# A COMPARISON OF CONTINUOUSLY CONTROLLED AND CONTROLLED K-THEORY

DOUGLAS R. ANDERSON
FRANCIS X. CONNOLLY
HANS J. MUNKHOLM

ABSTRACT. We define an unreduced version of the $\epsilon$-controlled lower $K$-theoretic groups of Ranicki and Yamasaki, [?], and Quinn, [?]. We show that the reduced versions of our groups coincide (in the inverse limit and its first derived, $\lim^1$) with those of [?]. We also relate the controlled groups to the continuously controlled groups of [?], and to the Quinn homology groups of [?].

Douglas R. Anderson was partially supported by the NSF (USA) under grant number DMS-9202598; Hans J. Munkholm was partially supported by the SNF (Denmark) under grants number 11-7792 and number 11-0062-1 PD.

<div align="center">CONTENTS</div>

<div align="center">1. INTRODUCTION.</div>

**1.1. The main results.** The main results are the following two theorems which appear (the first one in a slightly more general form) in Sections **??** and **??**, respectively. In them, the natural numbers do *not* include 0, i.e., $\mathbb{N} = \{1, 2, \cdots\}$.

**Theorem 1.1.** *Let $p : H \to B$ be a map from a topological space $H$ into a compact metric space $B$, and let $j \in \mathbb{N}$. There is a natural, short exact sequence*

$$0 \to \lim{}^1 K_{2-j}(B; p)_\epsilon \to K_{2-j}^{cc}(c(B_+), B_+; \overset{\circ}{c}(p_+)) \to K_{1-j}(B; p)_c \to 0.$$

*Assume further that each path component in $H$ and in $B$ is open. Then there is a reduced version of this seqeunce, i.e., a natural, short exact sequence*

$$0 \to \lim{}^1 \tilde{K}_{2-j}(B; p)_\epsilon \to \tilde{K}_{2-j}^{cc}(c(B_+), B_+; \overset{\circ}{c}(p_+)) \to \tilde{K}_{1-j}(B; p)_c \to 0.$$

**Theorem 1.2.** *Let $(\overline{Y}, C)$ be a pair of compact, metrizable spaces with $C$ tame in $\overline{Y}$. Let*

$$Y = \overline{Y} - C \overset{q}{\longleftarrow} Holink(C, \overline{Y}) \overset{p}{\longrightarrow} C$$

*be the correponding holink diagram. There are natural chain complexes of abelian groups*

$$K_2(Y) \overset{i_2}{\longrightarrow} K_2^{cc}(\overline{Y}, C; id_Y) \overset{\Delta_2}{\longrightarrow} K_1(C; p)_c \overset{\ell_1}{\longrightarrow} K_1(Y) \overset{i_1}{\longrightarrow} K_1^{cc}(\overline{Y}, C; id_Y) \overset{\Delta_1}{\longrightarrow} \cdots$$

*and*

$$Wh(C; p)_c \overset{\ell_1}{\longrightarrow} Wh(Y) \overset{i_1}{\longrightarrow} Wh^{cc}(\overline{Y}, C; id_Y) \overset{\Delta_1}{\longrightarrow} \tilde{K}_0(C; p)_c \overset{\ell_0}{\longrightarrow} \tilde{K}_0(Y) \overset{i_0}{\longrightarrow} \cdots.$$

*For each $j \in \mathbb{N}$, these complexes are exact at $K_{1-j}(C; p)_c$ and $K_{1-j}(Y)$, respectively at $\tilde{K}_{1-j}(Y; p)_c$ and $\tilde{K}_{1-j}(Y)$; the homology groups at $K_{2-j}^{cc}(\overline{Y}, C; id_Y)$ and at $\tilde{K}_{2-j}^{cc}(\overline{Y}, C; id_Y)$ are naturally isomorphic to $\lim^1 K_{2-j}(C; p)_\epsilon$.*

*For the reduced version there is the added hypothesis that each path component in each space in the holink diagram must be open.*

**Remarks: 1)** In Section **??**, we recall and comment upon a result of Ferry and Pedersen concerning the vanishing of the $lim^1$ groups above.
**2)** We strongly believe that the theorems hold also for $j = 0$, but we emphasize that the present paper does not prove this to be true.
**3)** The extra conditions on path components assumed for the reduced version are discussed in Section **??**.                                                        □

Let $p : E \to B = |K|$ be a map which admits an iterated mapping cylinder structure (also called a *homotopy colimit structure*, see [**?**]) over the finite simplicial

complex $K$, and let $S : \mathcal{TOP} \to \mathcal{SPEC}$ be the functor which associates to a space $X$ its non-connective $K$-theory spectrum as defined, e.g., in [**?**]. By Theorem 0.1 of [**?**], it is known that Quinn's homology groups $H_*(\Sigma K, v_+; \mathcal{S}(\Sigma p))$, derived from this functor $S$ as in the appendix of [**?**], are isomorphic to the boundedly controlled groups $K_*(p)$ defined in [**?**]. On the other hand, in [**?**], the latter groups are shown to coincide with the $cc$ $K$-theory groups occurring in Theorem 1.1. Thus, we have the following corollary.

**Corollary 1.3.** *For each $p$ as above, and each $j \in \mathbb{N}$, there is a short exact sequence*

$$0 \to \lim{}^1 K_{2-j}(|K|; p) \to H_{2-j}(\Sigma|K|, v_+; \mathcal{S}(\Sigma p)) \to K_{1-j}(|K|; p)_c \to 0.$$

The groups that appear in the above theorems can briefly be described as follows. For each $\epsilon > 0$, and each map $p : H \to B$, $\epsilon$-*controlled, reduced $K$-groups* $Wh(B; p)_\epsilon$ $(= \tilde{K}_1(B; p)_\epsilon)$ and $\tilde{K}_0(B; p)_\epsilon$ have been defined in a geometric way by Chapman, [**?**], and in an algebraic one by Quinn, [**?**], as well as [**?**]. Our $\epsilon$ decorated groups are unreduced versions of the groups defined in [**?**] and [**?**]. They form inverse systems under connecting homomorphisms $K_{2-j}(B; p)_\epsilon \to K_{2-j}(B; p)_\delta$ defined for $\epsilon \leq \delta$, and $K_{2-j}(B; p)_c$ is the inverse limit group $\lim K_{2-j}(B; p)_\epsilon$ (and similarly for reduced groups). Our definition is given in Section **??** in terms of homotopy classes of path matrix morphisms between geometric modules in the space $H$ and with "$\epsilon$-control on maps and homotopies" in the space $B$ via the map $p$.

The decoration $^{cc}$ means *continuously controlled*. Let $(\overline{Y}, C)$ be a pair as in Theorem 1.2, and let $R$ be any (reasonable) ring (which we shall usually suppress from the notation). In [**?**], the authors study the $K$-theory, $K_*^{cc}(\overline{Y}, C; R)$, of geometric $R$-modules (or, more generally, objects in any additive category) in the space $Y = \overline{Y} - C$ with morphisms which are "continuously controlled" at $C$. Although the ideas from [**?**] have had important applications (see e.g. [**?**] and [**?**]), from a geometric viewpoint it seems more natural to allow (homotopy classes of) path matrix morphisms between geometric modules in a space $H$ which comes with a map $p : H \to Y = \overline{Y} - C$. This map is then used to introduce control at $C$. This viewpoint is introduced in the continuously controlled setting by two of the present authors in [**?**].

In Theorem 1.1, the subscript $+$ refers to the addition of a disjoint base point. Also, $\mathring{c}(p_+)$ is the induced map from the open cone $\mathring{c}(H \amalg X)$ into $\mathring{c}(B_+) = c(B_+) - B_+$. Thus, the $cc$ groups appearing there are derived from the additive category of geometric modules in the open cone of $H \amalg X$ with morphisms that are $cc$ homotopy classes of $cc$ path matrix morphisms.

In Theorem 1.2, $id_Y$ is the identity map on $Y = \overline{Y} - C$, so the $cc$ groups appearing there date back to geometric modules in $Y$ with path matrix morphisms which are $cc$ at $C$ (up to homotopies which are also $cc$ at $C$).

We proceed to describe the background for our interest in these results.

**1.2. Background.** Both Quinn and Chapman were motivated by the controlled boundary problem for a manifold $M$, equipped with a map $p : M \to B$ with $B$ as above. The problem is to decide under which conditions $M$ admits a boundary in such a way that $p$ extends over $M \cup \partial M$. Under suitable a priori conditions on the dimension and the behaviour of $M$ near $\infty$, Chapman shows that the problem is completely governed by a primary obstruction and a secondary obstruction which lie in groups that are constructed via inverse limits and $\lim^1$ from suitable groups of the above $\epsilon$-controlled type. Quinn, in [?], also gave an obstruction theory solution to the controlled boundary problem, but the obstruction group is constructed (in Sections 5 and 8 in [?]) in a less transparent way. In his Notre Dame CBMS-lectures in 1984, Quinn outlined an algebraic definition of groups $Wh(B;p)_\epsilon$ and $\tilde{K}_0(B;p)_\epsilon$ in terms of geometric modules in $M$ and ($\epsilon$-small homotopy classes of) path matrix morphisms in $M$ with $\epsilon$-short image paths in $B$. These groups can be directly related to the cellular chains on $M$ (assuming a "controlled" cellular structure), and this introductory part of Quinn's lecture series was apparently geared towards a direct description (in terms of chain complexes) of the obstructions obtained in [?]. Unfortunately, this approach was never carried to fruition, although part of the projected CBMS volume appeared later as [?]. In fact, to the best of the present authors' knowledge,

(A) no one has ever established a strict link between the obstruction groups that appear in [?] and the groups defined in [?].

In particular, even though Quinn's obstruction groups do form a homology theory (as described in Section 8 of [?]),

(B) there is no documentation in the literature to support the folklore belief that in some great generality the inverse limit groups $Wh(B;p)_c$, $\tilde{K}_0(B;p)_c$, $\cdots$ are the groups of a homology theory.

Also,

(C) no solid foundation for an interpretation of Quinn's obstructions in terms of chain complexes in some category of geometric modules seems to exist.

Along a similar vein, apparently

(D) nobody has taken time out to construct the expected isomorphism between Chapman's geometrically defined and Quinn's algebraically defined $\epsilon$-controlled groups, and to study its properties.

The present paper is related to all of these problems. Starting from behind, (D) should actually be very accessible now that the algebraic version of the groups involved is equally well understood (thanks to [?] and the present paper) as is the original, geometric version in [?].

Since Chapman's obstructions are quite explicitly (and classically) given in [?], the precise form of a solution to (D) might actually be of value in trying to understand (C). Since the primary obstruction is in a group of type $\tilde{K}_0(B;p)_c$ and the secondary in one of type $\lim^1 Wh(B; p)_\epsilon$, Theorem 1.1 suggests a total obstruction in a certain $cc$ Whitehead group $Wh^{cc}(\overline{M}, B; id_M)$; see Remark **4)** following **??**.

Recently, in [**?**], Ranicki and Yamasaki have given a chain complex version of Quinn's $\epsilon$-controlled groups from [**?**]. Their main goal is to give an algebraic proof of the result of Chapman [**?**], that Whitehead torsion is a topological invariant. As such it represents the culmination of the program initiated by Connell and Hollingsworth, in [**?**], with the invention of the notion of a geometric group. They also reprove other results, including the Borsuk conjecture (originally proved by West, in [**?**]). Their work includes "stably exact" Mayer-Vietoris sequences which are, of course, very much related to our problem (B), but they do not attack (B) directly. In [**?**], it is shown that the $cc$ groups appearing in Theorem 1.1 do form a homology theory, defined on a category of maps like $p$ (or, more generally, on a category of diagrams of holink type). Thus, the controlled groups $K_{2-j}(B; p)_c$ have homology theory behavior, provided the $\lim^1$-terms in Theorem 1.1 vanish. We discuss a vanishing result (due essentially to Ferry and Pedersen) briefly in Section **??**.

### 1.3. Historical remarks and acknowledgements.

The exact sequences in Theorem 1.1 were conjectured to be true (with a *boundedly controlled* group instead of the $cc$ group) a long time ago by two of the authors, DRA and HJM. Actually, HJM's lecture at the Oberwolfach Topologie Tagung in 1987, contained an explicit version of this conjecture. One ingredient in the proof (injectivity of $\lim^1 K_1(B; p)_\epsilon \to K_1^{cc}(c(B_+), B_+; \overset{\circ}{c}(p_+))$) eluded them and the project was essentially shelved until HJM came to University of Notre Dame on a sabbatical in 1992/93. There, discussions between him and the third author, FXC, about the groups $K^{cc}(\overline{Y}, C; id_Y)$ appearing in Theorem 1.2, renewed the interest and the missing part of the proof was finally found.

DRA wants to thank the Mathematics Department at the University of Notre Dame and the Mathematical Sciences Research Institute, Berkely, for their kind hospitality during different periods when this paper was in preparation.

HJM wants to thank the Mathematics Department of University of Notre Dame for a pleasant atmosphere (mathematical and otherwise) during the above mentioned sabbatical. Special thanks go to Bruce Williams for his help in arranging the sabbatical and for his never ending willingness to discuss mathematics and to search through his extensive files of (p)reprints to find the one you want, or even the one you *ought* to want. He also thanks the Mathematical Sciences Research Institute, Berkeley, for its hospitality during the summer of 1993 when much of the actual writing was done. The delay in the appearance of this preprint arises from the authors' desire to have [**?**] in final form before using results from it.

## 2. Geometric modules and path matrix morphisms.

In this section we set up the general framework (a priori involving no control) for *geometric modules, path matrix morphisms, and path matrix homotopies* in any topological space $E$.

We recall that a *Moore path* in $E$ is a pair $(\omega, l_\omega)$ where $\omega : [0, \infty) \to E$ is continuous, and constant on the subray $[l_\omega, \infty)$. A homotopy between Moore paths is a pair $(H, l_H)$ where $l_H : [0, 1] \to [0, \infty)$ and $H : [0, \infty) \times [0, 1] \to E$ are continuous, and $H(-, y)$ is constant on $[l_H(y), \infty)$ for each $y \in [0, 1]$. As usual, $(H, l_H)$ is a homotopy *from* $\delta_1(H, l_H) = (H(-, 0), l_H(0))$ *to* $\delta_0(H, l_H) = (H(-, 1), l_H(1))$. Often $l_\omega$ and $l_H$ will be suppressed from the notation, $\omega(l_\omega)$ will be denoted $\omega(\infty)$, and $H(l_H(y), y)$ will be denoted $H(\infty, y)$. The set, $\mathcal{M}(E)$, of Moore paths in $E$ will be viewed as a category whose objects are the points of $E$ and where the morphisms from $x$ to $y$ are the Moore paths $\omega$ with $y = \omega(0)$ and $x = \omega(\infty)$ (so that the usual concatenation of paths defines the composition). Also, $\mathcal{H}(E)$ denotes the set of end point fixing homotopies, so $\delta_i$ is a map $\mathcal{H}(E) \to \mathcal{M}(E)$, for $i = 0, 1$.

A *geometric module* on $E$ is a pair $(S, \sigma)$ with S a discrete set and $\sigma : S \to E$ a function. We think of each $s \in S$ as a basis element which "sits at" the point $\sigma(s) \in E$. A *(path matrix) morphism* $\varphi : (T, \tau) \to (S, \sigma)$ is a function[1] $\varphi : S \times T \to R\mathcal{M}(E)$ into the free $R$-module generated by $\mathcal{M}(E)$. There are two requirements on $\varphi(s, t) = \Sigma_\omega r_{st}(\omega)\omega$, viz.

(1) $\qquad\qquad r_{st}(\omega) \neq 0 \Rightarrow [\omega(0) = \sigma(s) \text{ and } \omega(\infty) = \tau(t)],$

and

(2) $\qquad\qquad \forall t \in T : \{(s, \omega) \in S \times \mathcal{M}(E) \,|\, r_{st}(\omega) \neq 0\}$ is finite.

The former condition says that each Moore path is to be considered a morphism from its end point to its beginning point; general morphisms are then $R$-linear combinations of such "atomic" morphisms. The second condition guarantees that matrix multiplication and the composition in $\mathcal{M}(E)$ give a well defined composition of morphisms. The resulting category will be denoted $GM(E)$. It is additive, with addition of morphisms coming from the abelian group structure in $R\mathcal{M}(E)$ and with $(S_1, \sigma_1) \oplus (S_2, \sigma_2)$ given by the disjoint union $(S_1 \coprod S_2, \sigma_1 \coprod \sigma_2)$. The inclusions and projections

$$(S_i, \sigma_i) \xrightarrow{\iota_i} (S_1, \sigma_1) \oplus (S_2, \sigma_2) \xrightarrow{\pi_i} (S_i, \sigma_i), \, i = 1, 2,$$

involve only constant paths $\gamma_x$, $x \in E$. For example, $\iota_2(s, s_2) = \delta_{s,s_2}\gamma_{\sigma_2(s_2)}$ for $(s, s_2) \in (S_1 \coprod S_2) \times S_2$, where $\delta$ is the Kronecker delta. Incidentally, if $S_1 = \emptyset$, this defines the identity morphism on $(S_2, \sigma_2)$.

---

[1] The indexing corresponds to classical matrix notation, i.e., the row index $s$ corresponds to the range, and the column index $t$ to the domain.

We call a morphism $\varphi : M = (T, \tau) \to N = (S, \sigma)$ *geometric*[2] if it has the form

$$(3) \qquad\qquad \varphi(t, s) = \delta_{t, \alpha(s)} u_s \omega_s,$$

where $\alpha : S \to T$ is a map; $u_s \in R$; and $\omega_s \in \mathcal{M}(E)$ has $\omega_s(\infty) = \tau\alpha(s), \omega_s(0) = \sigma(s)$[3]. We note that the above inclusions $\iota_i$ and projections $\pi_i$ are examples of such geometric morphisms.

If $\alpha$ is a bijection, each $u_s$ is a unit in $R$, and each $\omega_s$ is constant, then the corresponding $\varphi$ is an isomorphism in $GM(E)$ with an inverse obtained by inverting $\alpha$ and each $u_s$. Clearly, the direct sum operation $\oplus$ is associative and commutative up to coherent, geometric isomorphisms of this sort (with all $u_* = 1$). By abuse of set theory, we shall actually proceed as if $\oplus$ is strictly associative, i.e., as if $\coprod$ is strictly associative in the category of sets.

A *(path matrix) homotopy* from $(T, \tau)$ to $(S, \sigma)$ is a function $\Phi : S \times T \to R\mathcal{H}(E)$, say $\Phi(s, t) = \Sigma_H r_{st}(H)H$, which satisfies

$$(4) \qquad r_{st}(H) \neq 0 \Rightarrow [H(0, y) = \sigma(s) \text{ and } H(\infty, y) = \tau(t), \forall y \in [0, 1]],$$

and

$$(5) \qquad\qquad \forall t \in T : \{(s, H) \in S \times \mathcal{H}(E) \,|\, r_{st}(H) \neq 0\} \text{ is finite.}$$

The end point maps $\delta_i : \mathcal{H}(E) \to \mathcal{M}(E)$, $i = 0, 1$, extend by linearity to maps $\delta_i : R\mathcal{H}(E) \to R\mathcal{M}(E)$, and $\Phi$ is understood to be a homotopy from $\delta_1\Phi$ to $\delta_0\Phi$. We also write $\Phi : \delta_1\Phi \simeq \delta_0\Phi$. This homotopy notion behaves well viz-a-viz composition, sum and direct sums so that one has a corresponding, additive, homotopy category $\mathcal{GM}(E)$.

In this generality, $GM(E)$ and $\mathcal{GM}(E)$ are of little interest. To arrive at interesting constructions, one has to impose conditions on the objects, the morphisms, and the homotopies allowed: For objects one has some "finiteness" condition on $\sigma : S \to E$. For morphisms $\varphi$ and homotopies $\Phi$ the conditions are expressed in terms of the following families of subsets of $E$ (where $\varphi(s, t) = \Sigma r_{st}(\omega)\omega$ and $\Phi(s, t) = \Sigma r_{st}(H)H$).

$$(6) \qquad \mathcal{F}(\varphi) \;=\; \{\omega([0, \infty)) \subseteq E \,|\, \exists (s, t) \in S \times T : r_{st}(\omega) \neq 0\},$$

$$(7) \qquad \mathcal{F}(\Phi) \;=\; \{H([0, \infty) \times [0, 1]) \subseteq E \,|\, \exists (s, t) \in S \times T : r_{st}(H) \neq 0\}.$$

In Section **??** we study one such set of conditions ("$\epsilon$-control") in some detail. Another set of conditions ("continuous control") is studied by two of the present authors in [?], in part in order to obtain the *cc* results needed in the present paper. Some definitions and results from [?] are summarized in Section **??**.

---

[2]In [**?**], such a morphism (or rather its homotopy class) is said to be "induced by a basis morphism". We have borrowed the present term from [**?**] and [**?**]. However, we use the term "geometric" in connection with morphisms only in the context defined here, in contradistinction to [**?**], where all path matrix morphisms are called "geometric morphisms"; but where "geometric" is also used in our, restricted sense.

[3]i.e., for given $s$, $\sigma(s)$ is the beginning point of at most one path having non zero coefficient in $\varphi$; this path ends at $\tau\alpha(s)$ and has multiplicity $u_s$.

## 3. Controlled $K$-theory.

The category $\mathcal{TOP}_c$ of *controlled spaces* has as its objects all pairs $(B; p)$ where $p : E \to B$ be a map from an arbitrary space $E$ into a compact metric space $B$. A morphism $(f, g) : (B; p) \to (B'; p')$ is defined to be a pair of continuous maps

$$(8) \qquad (f : E \to E', g : B \to B') \text{ with } p'f = gp.$$

Let $j \in \mathbb{N}$ and $\epsilon > 0$. In this section we define and study the $\epsilon$-*controlled groups* $K_{2-j}(B; p)_\epsilon$, (which are *not* functorial on $\mathcal{TOP}_c$) as well as the *controlled groups* $K_{2-j}(B; p)_c$ ($= \lim K_{2-j}(B; p)_\epsilon$), and the *lim$^1$-groups* $\lim^1 K_{2-j}(B; p)_\epsilon$ (which *are* functorial on $\mathcal{TOP}_c$).

The first subsection deals with the cases $j = 1$ and $j = 2$, i.e., with $K_1-$ and $K_0-$groups. Reduced versions of these are defined in Section **??**, and compared to the similar groups defined by Ranicki and Yamasaki ([?]) in Section **??**. Finally, in Section **??**, results from [?] and Section **??** are utilized to define the *lower* groups $K_{2-j}(B; p)_\epsilon$, $K_{2-j}(B; p)_c$ and $\lim^1 K_{2-j}(B; p)_\epsilon$ for $j > 2$. This uses a geometric form of Bass's contracted functor viewpoint, cf. [?].

First, however, we must discuss the relevant conditions on modules, morphisms and homotopies. All geometric modules $\sigma : S \to E$ will have $S$ *finite*. They will often be denoted by $M, N, L, \cdots$ instead of $(S, \sigma), (T, \tau), (R, \rho), \cdots$. The families $\mathcal{F}(\varphi)$ and $\mathcal{F}(\Phi)$ associated to a morphism $\varphi$ or a homotopy $\Phi$ in $GM(E)$ (see (**??**) and (**??**)) are *end pointed* in the sense that each $F \in \mathcal{F}$ ($= \mathcal{F}(\varphi)$ or $\mathcal{F}(\Phi)$) comes equipped with an ordered pair of points $(b_F, e_F)$. Such an end pointed family is said to be *bounded by $\epsilon \geq 0$*, if

$$(9) \qquad \forall F \in \mathcal{F} : p(F) \subseteq B(p(b_F), \epsilon) \cap B(p(e_F), \epsilon),$$

where $B(x, \epsilon)$ denotes the closed ball of radius $\epsilon$ and with center $x \in B$. An $\epsilon$-*morphism* $\varphi : N \to M$ between geometric modules on $E$ is a morphism for which $\mathcal{F}(\varphi)$ is bounded by $\epsilon$. Similarly, one defines an $\epsilon$-homotopy between $\epsilon$-morphisms, written $\varphi \simeq_\epsilon \psi$. The collection $GM_\epsilon(B; p)$ of all (finite) geometric modules on $E$ and all $\epsilon$-morphisms between them, is not a subcategory of $GM(E)$. The best one can say is that the composition maps $GM_\delta(B; p)(M, L) \times GM_\epsilon(B; p)(N, M)$ into $GM_{\delta+\epsilon}(B; p)(N, L)$, and it is compatible with the inclusions $GM_\epsilon(B; p) \subseteq GM_{\epsilon'}(B; p)$ for $\epsilon \leq \epsilon'$. On the other hand, $\epsilon$-homotopy *is* an equivalence relation in the collection $GM_\epsilon(B; p)$, so one does have the corresponding quotient "non-categories" $\mathcal{GM}_\epsilon(B; p)$ with compositions

$$\mathcal{GM}_\delta(B; p)(M, L) \times \mathcal{GM}_\epsilon(B; p)(N, M) \to \mathcal{GM}_{\delta+\epsilon}(B; p)(N, L).$$

Also, $\epsilon$-homotopies can be composed with $\delta$-morphisms or $\delta$-homotopies on either side to give $(\epsilon + \delta)$-homotopies. We leave it to the reader to keep track of the obvious formal properties of such compositions, cf. Section 2 of [**?**].

An $\epsilon$-morphism $\xi : N \to N$ is called an $\epsilon$-*idempotent*, if $\xi^2 \simeq_{2\epsilon} \xi$. If $\xi$ is an $\epsilon$-idempotent, then so is $1_N - \xi$. A *geometric $\epsilon$-idempotent* is an $\epsilon$-morphism which is geometric in the sense of (**??**), and for which $\alpha = 1_T : S = T \to T$ and each $u_s \in \{0, 1\}$.

An $\epsilon$-morphism $\varphi : N \to M$ is called an $\epsilon$-*isomorphism*, if there exists an $\epsilon$-morphism $\varphi' : M \to N$ such that $\varphi'\varphi \simeq_{2\epsilon} 1_N$ and $\varphi\varphi' \simeq_{2\epsilon} 1_M$. We also call $\varphi'$ *an $\epsilon$-inverse* for $\varphi$. We note that an $\epsilon$-morphism, $\varphi$, which is geometric and for which the map $\alpha : S \to T$ is a bijection and each $u_s$ is a unit in $R$, is an $\epsilon$-isomorphism, with a preferred, geometric $\epsilon$-inverse, $\varphi^{-1}$, obtained by inverting $\alpha$, reversing all the paths involved, and inverting all the units $u_s$. We shall refer to such a morphism as a *geometric $\epsilon$-isomorphism.*

**Remark:** From a purely algebraic viewpoint, it might be natural to allow each $u_s$ to be an arbitrary idempotent in $R$ (respectively, an arbitrary unit in $R$) in the above definition of a geometric idempotent (respectively, a geometric isomorphism). However, in geometry the above choices seem preferable, cf. Section **??**.  $\square$

Clearly, any geometric module $M = (S, \sigma)$ admits a geometric 0-isomorphism to a module $(T, \tau)$ with $T \subset E \times \mathbb{N}$ and $\tau = proj_1 | T : T \to E$. In fact, one may take each path involved in such an isomorphism to be constant. Therefore, up to conjugation by geometric 0-isomorphisms, there is only a *set* of $\epsilon$-automorphisms and a *set* of $\epsilon$-idempotents in $GM_\epsilon$. This is important for set theoretic reasons in the following section.

**3.1. Definition of $K_i(B; p)_\epsilon$ for $i = 0, 1$.** The group $K_0(B; p)_\epsilon$ will be defined in terms of the class of $\epsilon$-idempotents in $GM_\epsilon(B; p)$. Thus, let $\xi : M \to M$ and $\eta : N \to N$ be $\epsilon$-idempotents. We write $\xi \sim_\epsilon \eta$, if there exist a geometric module $L$ and $\epsilon$-morphisms $\varphi : M \oplus L \to N \oplus L$, $\psi : N \oplus L \to M \oplus L$ such that

$$(10) \qquad \psi\varphi \simeq_{2\epsilon} \xi \oplus 1_L \text{ and } \varphi\psi \simeq_{2\epsilon} \eta \oplus 1_L;$$

and

$$(11) \qquad \psi \simeq_{2\epsilon} (\xi \oplus 1_L)\psi \simeq_{2\epsilon} \psi(\eta \oplus 1_L) \text{ and } \varphi \simeq_{2\epsilon} (\eta \oplus 1_L)\varphi \simeq_{2\epsilon} \varphi(\xi \oplus 1_L).$$

The idea here is that up to homotopy $\varphi$ is an isomorphism $(M \oplus L, \xi \oplus 1_L) \to (N \oplus L, \eta \oplus 1_L)$ between objects in the idempotent completion of $GM_\epsilon(B; p)$ and that $\psi$ is its inverse. We note that $\sim_\epsilon$ is not an equivalence relation, but it is symmetric and reflexive, and it does have the following properties

$$(12) \qquad [\xi \sim_\epsilon \eta \text{ and } \eta \sim_\epsilon \zeta] \Rightarrow \xi \sim_{5\epsilon} \zeta,$$

(13)                     $\forall \epsilon$-idempotent $\xi : M \to M : \xi \oplus (1_M - \xi) \sim_\epsilon 1_M,$

(14)                     $[\xi_i \sim_\epsilon \eta_i, i = 1, 2] \Rightarrow [\xi_1 \oplus \xi_2 \sim_\epsilon \eta_1 \oplus \eta_2],$

and

(15)                     $[\xi \sim_\delta \eta \text{ and } \delta < \epsilon] \Rightarrow [\xi \sim_\epsilon \eta].$

Let $\rightsquigarrow_\epsilon$ denote the equivalence generated by $\sim_\epsilon$. Since $\epsilon$-idempotents which are conjugate by means of a geometric 0-isomorphism, are $\sim_\epsilon$-equivalent, the $\rightsquigarrow_\epsilon$-equivalence classes form a set. Since $\oplus$ is associative and commutative up to geometric 0-isomorphisms, this set is an abelian semigroup under a sum induced by $\oplus$. We define $K_0(B; p)_\epsilon$ to be the corresponding abelian group. Thus, any $x \in K_0(B; p)_\epsilon$ has the form $x = [\xi] - [\eta]$ with $\xi$ and $\eta$ $\epsilon$-idempotents, and $[\xi] - [\eta] = [\xi'] - [\eta']$ if and only if there exist $\epsilon$-idempotents $\xi_1, \xi_2, \cdots, \xi_{m+1}$, and $\zeta$, such that

(16)        $\xi \oplus \eta' \oplus \zeta = \xi_1 \sim_\epsilon \xi_2 \sim_\epsilon \xi_3 \sim_\epsilon \cdots \sim_\epsilon \xi_{m+1} = \xi' \oplus \eta \oplus \zeta.$

Since there is no bound on $m$, it seems that there is no control on the relation between $\xi \oplus \eta'$ and $\xi' \oplus \eta$. However, by using **(??)**–**(??)**, one can imitate Theorem 3.5 of [**?**] (see also Lemma 4.4 of [**?**]) to get the following equivalences (where $\bar{\xi}$ is short for $1 - \xi$)

$$\xi_1 \oplus 1 \oplus \cdots \oplus 1 \oplus 1 \sim_\epsilon$$
$$\xi_1 \oplus \bar{\xi}_2 \oplus \xi_2 \oplus \cdots \oplus \bar{\xi}_{m-1} \oplus \xi_{m-1} \oplus \bar{\xi}_m \oplus \xi_m \sim_\epsilon$$
$$\xi_2 \oplus \bar{\xi}_2 \oplus \xi_3 \oplus \cdots \oplus \bar{\xi}_{m-1} \oplus \xi_m \oplus \bar{\xi}_m \oplus \xi_{m+1} \sim_\epsilon$$
$$1 \oplus 1 \oplus \cdots \oplus 1 \oplus 1 \oplus \xi_{m+1} \sim_0 \xi_{m+1} \oplus 1 \oplus 1 \oplus \cdots \oplus 1 \oplus 1.$$

In view of **(??)** and **(??)**, this proves the following proposition.

**Proposition 3.1.** *There is a universal constant $k \in \mathbb{N}$ such that whenever $\xi$, $\eta$, $\xi'$, and $\eta'$ are $\epsilon$-idempotents with $[\xi] - [\eta] = [\xi'] - [\eta'] \in K_0(B; p)_\epsilon$, then $\xi \oplus \eta' \sim_{k\epsilon} \xi' \oplus \eta'$.*

The inclusions $GM_\delta \subseteq GM_\epsilon$, $\delta < \epsilon$, give rise to "relax control" homomorphisms $K_0(B; p)_\delta \to K_0(B; p)_\epsilon$. We define the *controlled group* $K_0(B; p)_c$ to be the (inverse) limit of the resulting inverse system indexed on all $\epsilon > 0$. We also need the first derived, $\lim^1 K_0(B; p)_\epsilon$.

We note that $K_0(B; p)_\epsilon$ is *not functorial* on the category $\mathcal{TOP}_c$. However, since $B$ is assumed compact, the inverse system $\{K_0(B; p)_\epsilon\}_\epsilon$, when viewed as an object in the usual pro category (cf. the Appendix in [**?**]), *is* functorial. Hence the controlled groups $K_0(B; p)_c$ and the groups $\lim^1 K_0(B; p)_\epsilon$ are also functorial.

The groups $K_1(B;p)_\epsilon$ and $K_1(B;p)_c$ are constructed in a similar way, starting from the class of $\epsilon$-automorphisms in $GM_\epsilon(B;p)$. An $\epsilon$-automorphism $\gamma : N \to N$ is called $\epsilon$-elementary, if $N$ admits a decomposition[4] $N = N_1 \oplus N_2$ so that $\gamma$ is of the form

$$(17) \qquad \gamma = \begin{pmatrix} 1 & \gamma' \\ 0 & 1 \end{pmatrix} : N_1 \oplus N_2 \to N_1 \oplus N_2$$

for some $\epsilon$-morphism $\gamma' : N_2 \to N_1$. Such a $\gamma$ has a preferred $\epsilon$-inverse

$$(18) \qquad \gamma^{-1} = \begin{pmatrix} 1 & -\gamma' \\ 0 & 1 \end{pmatrix} : N_1 \oplus N_2 \to N_1 \oplus N_2.$$

which is itself $\epsilon$-elementary and has $\gamma\gamma^{-1} = \gamma^{-1}\gamma = 1$. Let $\varphi : M \to M$ and $\psi : N \to N$ be $\epsilon$-automorphisms. We write $\varphi \sim_\epsilon \psi$, if there exist geometric modules $L', L''$; a geometric 0-isomorphism $\nu : M \oplus L' \to N \oplus L''$; an $\epsilon$-elementary automorphism $\gamma : N \oplus L'' \to N \oplus L''$; and a $2\epsilon$-homotopy

$$(19) \qquad \nu(\varphi \oplus 1_{L'})\nu^{-1} \simeq_{2\epsilon} (\psi \oplus 1_{L''})\gamma.$$

The relation thus defined is not an equivalence relation, but it does satisfy

$$(20) \qquad \varphi \sim_\epsilon \psi \Rightarrow \psi \sim_{3\epsilon/2} \varphi,$$

$$(21) \qquad [\varphi_i \sim_\epsilon \psi_i, i = 1, 2] \Rightarrow \varphi_1 \oplus \varphi_2 \sim_\epsilon \psi_1 \oplus \psi_2,$$

and

$$(22) \qquad [\varphi \sim_\delta \psi \text{ and } \delta < \epsilon] \Rightarrow \varphi \sim_\epsilon \psi.$$

Let $\rightsquigarrow_\epsilon$ be the equivalence relation generated by $\sim_\epsilon$. The equivalence classes form (a set, and) an abelian semigroup under a sum induced by $\oplus$. The zero element is the class of $1_M$ for any $M$.

If $\varphi$ has $\varphi'$ as an $\epsilon$-inverse, then successive multiplication on the right by the following six $\epsilon$-elementary automorphisms

$$(23) \qquad \begin{pmatrix} 1 & 0 \\ 1 & 1 \end{pmatrix}, \begin{pmatrix} 1 & -1 \\ 0 & 1 \end{pmatrix}, \begin{pmatrix} 1 & 0 \\ 1 & 1 \end{pmatrix}, \begin{pmatrix} 1 & 0 \\ -\varphi & 1 \end{pmatrix}, \begin{pmatrix} 1 & \varphi' \\ 0 & 1 \end{pmatrix}, \begin{pmatrix} 1 & 0 \\ -\varphi & 1 \end{pmatrix}$$

gives a string of six $\sim_\epsilon$-equivalences

$$(24) \qquad \begin{pmatrix} 1 & 0 \\ 0 & 1 \end{pmatrix} \sim_\epsilon \cdots \sim_\epsilon \begin{pmatrix} \varphi & 0 \\ 0 & \varphi' \end{pmatrix}.$$

This establishes the existence of inverses, and the group so defined will be called $K_1(B;p)_\epsilon$.

---

[4]Recall that if $N = (T, \tau)$, then this means that $N_i = (T_i, \tau|T_i)$ where $T = T_1 \coprod T_2$.

The following proposition is essentially a version of Lemma 4.4 of [?] or Theorem 3.5 of [?].

**Proposition 3.2.** *There is a universal constant $k \in \mathbb{N}$ such that whenever two $\epsilon$-automorphisms $\varphi : M \to M$ and $\psi : N \to N$ represent the same element in $K_1(B; p)_\epsilon$, then there exist a geometric module $L$ and thirteen $\epsilon$-elementary automorphisms $\Gamma^i : M \oplus L \oplus N \to M \oplus L \oplus N$, $i = 1, 2, \cdots, 13$, such that*

$$\varphi \oplus 1_L \oplus 1_N \simeq_{k\epsilon} (1_M \oplus 1_L \oplus \psi) \Gamma^1 \Gamma^2 \cdots \Gamma^{13}.$$

**Proof:** Since conjugation by geometric 0-isomorphisms preserves $\epsilon$-elementary automorphisms, we do not have to worry about associativity and commutativity isomorphisms in the proof. Using (**??**), we see that there is an integer $m$ and a sequence of $3\epsilon$-homotopies $\nu_i(\varphi_i \oplus 1_{L'_i})\nu_i^{-1} \simeq_{3\epsilon} (\varphi_{i+1} \oplus 1_{L''_i})\gamma^i, i = 1, 2, \cdots, m-1$, where each $\varphi_i$ is an $\epsilon$-automorphism, each $\nu_i$ a geometric 0-isomorphism, and each $\gamma^i$ an $\epsilon$-elementary[5] automorphism. Moreover, $\varphi_1 = \varphi$ and $\varphi_m = \psi$. After replacing each $\varphi_i$ by $1_{L''_1 \oplus \cdots \oplus L''_{i-1}} \oplus \varphi_i \oplus 1_{L'_i \oplus \cdots \oplus L'_{m-1}}$ (and adjusting each $\nu_i$ and $\gamma^i$ correspondingly), we may assume that $L'_i = L''_i = 0$ for all $i$. This, of course, also changes $\varphi$ and $\psi$, but the conclusion for the new version does imply that for the old one. Next, we can replace each (new) $\varphi_i$ by $\nu_{[i]}\varphi_1\nu_{[i]}^{-1}$, where $\nu_{[i]} = \nu_{m-1}\nu_{m-2}\cdots\nu_i$, and each $\nu_i$ by 1. Again, this changes $\varphi$ and $\psi$, but the conclusion for the new pair implies that for the old. We now have $\varphi_i \simeq_{3\epsilon} \varphi_{i+1}\gamma^i, i = 1, 2, \cdots m-1$, from which we conclude that

$$(25) \qquad \varphi_1 \oplus [\oplus_{i=2}^{m-1}(\varphi'_i \oplus \varphi_i)] \simeq_{3\epsilon} [\oplus_{i=2}^{m-1}(\varphi_i \oplus \varphi'_i) \oplus \varphi_m]\delta,$$

where $\varphi'_i$ is an $\epsilon$-inverse of $\varphi_i$ and $\delta = \oplus_{i=1}^{m-2}(\gamma^i \oplus 1) \oplus \gamma^{m-1}$ is an $\epsilon$-elementary automorphism. For any $\epsilon$-morphism $\varphi$ with $\epsilon$-inverse $\varphi'$, the product of the six factors in (**??**) is $3\epsilon$-homotopic to $\varphi \oplus \varphi'$. Therefore, there are twelve $\epsilon$-elementary automorphisms $\delta^i$ such that

$$(26) \qquad 1_M \oplus [\oplus_{i=2}^{m-1}(\varphi'_i \oplus \varphi_i] \simeq_{3\epsilon} \prod_{i=1}^{6}\delta^i, \text{ and } [\oplus_{i=2}^{m-1}\varphi_i \oplus \varphi'_i] \oplus 1_N \simeq_{3\epsilon} \prod_{i=7}^{12}\delta^i.$$

The desired conclusion follows easily from (**??**) and (**??**). $\qquad \square$

If $\delta < \epsilon$, then by (**??**), the inclusions $GM_\delta \subseteq GM_\epsilon$ give rise to "relax control" homomorphisms $K_1(B; p)_\delta \to K_1(B; p)_\epsilon$. We define the *controlled group* $K_1(B; p)_c$ to be the (inverse) limit of the resulting inverse system indexed on all $\epsilon > 0$. We also need the first derived, $\lim^1 K_1(B; p)_\epsilon$.

As with $K_0$, the individual groups $K_1(B; p)_\epsilon$ are *not* functorial, but the inverse system $\{K_1(B; p)_\epsilon\}_\epsilon$ and the groups $K_1(B; p)_c$ and $\lim^1 K_1(B; p)_\epsilon$ *are* functorial on $\mathcal{TOP}_c$.

---

[5]The wording of (**??**) would give the bound $3\epsilon/2$ only, but the easy proof of (**??**) does give $\epsilon$ here.

**3.2. The reduced groups for** $i = 0, 1$. We recall the notion of a *geometric* morphism, defined in (**??**), and denote by $GK_i(B; p)_\epsilon$ the subgroup of $K_i(B; p)_\epsilon$ generated by all $[\varphi]$ where $\varphi$ is a geometric $\epsilon$-idempotent (if $i = 0$), or a geometric $\epsilon$-automorphism (if $i = 1$), cf. Section **??**. The reduced groups $\tilde{K}_i(B; p)_\epsilon$, are defined to fit into the short exact sequences.

$$(27) \qquad 0 \to GK_i(B; p)_\epsilon \to K_i(B; p)_\epsilon \to \tilde{K}_i(B; p)_\epsilon \to 0,$$

$i = 0, 1$. For $i = 1$, we also write $Wh(B; p)_\epsilon$ for the reduced group. In order to study the behavior of the reduced groups under lim and $\lim^1$, we need two lemmas.

**Lemma 3.3.** *If two geometric modules are $\epsilon$-isomorphic, then there is a geometric $\epsilon$-isomorphism between them.*

**Proof:** One can easily adapt the proof of Theorem 2.6, pp. 140ff. of [**?**] to the present situation. $\square$

**Lemma 3.4.** *For any geometric $\epsilon$-isomorphism $\nu : M \to N$ and any $\delta > 0$, there exist an $m \in \mathbb{N}$; geometric modules $M_1 = M, M_2, \cdots, M_m = N$; and geometric $\delta$-isomorphisms $\nu_i : M_i \to M_{i+1}, i = 1, 2, \cdots, m-1$, such that $\nu = \nu_{m-1} \cdots \nu_2 \nu_1$.*

**Proof:** One uses a straigthforward subdivision of paths (cf. Section 4 of [**?**]). It is, of course, essential that $B$ is compact. $\square$

**Proposition 3.5.** *The inverse system $\{GK_0(B; p)_\epsilon\}_\epsilon$ is stably isomorphic to the constant system $H_0(E, \mathbb{Z})$, where $H_0$ denotes singular homology. Consequently, one has the following natural isomorphisms and natural short exact sequence*

$$\lim^1 GK_0(B; p)_\epsilon = 0 \ , \ \lim^1 K_0(B; p)_\epsilon \cong \lim^1 \tilde{K}_0(B; p)_\epsilon,$$

$$0 \to H_0(E; \mathbb{Z}) \to K_0(B; p)_c \to \tilde{K}_0(B; p)_c \to 0.$$

*The inverse system $\{GK_1(B; p)_\epsilon\}_\epsilon$ satisfies the Mittag-Leffler condition. Consequently, one has the following natural isomorphisms and natural short exact sequence.*

$$\lim^1 GK_1(B; p)_\epsilon = 0 \ , \ \lim^1 K_1(B; p)_\epsilon \cong \lim^1 Wh(B; p)_\epsilon,$$

$$0 \to \lim GK_1(B; p)_\epsilon \to K_1(B; p)_c \to Wh(B; p)_c \to 0.$$

**Proof:** Any geometric $\epsilon$-idempotent $\xi : M \to M$ is $\sim_\epsilon$-equivalent to a geometric $\epsilon$-idempotent $\xi' : M' \to M'$ which is an $\epsilon$-isomorphism (i.e., each of the elements $u_s$ occurring in the presentation of $\xi'$ equals 1). In fact, if $M = (S, \sigma)$, one can take $\xi' = \xi|M'$ where $M' = (S', \sigma|S')$ and $S' = \{s \in S|u_s = 1\}$. Thus we need only consider geometric $\epsilon$-idempotents $\xi : M \to M$ which have $\epsilon$-inverses (denoted $\xi'$). Then $\xi \simeq_{3\epsilon} \xi\xi\xi' \simeq_{3\epsilon} \xi\xi' \simeq_{2\epsilon} 1$, so $\xi \simeq_{3\epsilon} 1_M$. Hence, if two $\epsilon$-invertible $\epsilon$-idempotents, $\xi : M \to M$ and $\eta : N \to N$, have $\xi \sim_\epsilon \eta$ so that (**??**) and (**??**)

hold, then $\psi$ and $\varphi$ are $\frac{3\epsilon}{2}$-inverses of one another. Then by Lemma **??**, there is a geometric $\frac{3\epsilon}{2}$-isomorphism $M \oplus L \to N \oplus L$. Therefore, there is a homomorphism $\rho_\epsilon : GK_0(B;p)_\epsilon \to H_0(E, \mathbb{Z})$ mapping $[\xi]$ to $\Sigma n_\beta [e_\beta]$ where $M = (S, \sigma)$, $[e_\beta]$ is the canonical, zero dimensional homology class in the path component $E_\beta \subseteq E$, and $n_\beta$ is the cardinality of $\sigma^{-1}(E_\beta)$. Conversely, if we fix the points $\{e_\beta\}$, there is an obvious homomorphism $\rho_\epsilon' : H_0(E; \mathbb{Z}) \to GK_0(B;p)_\epsilon$ with $\rho_\epsilon'([e_\beta]) = 1_{(T,\tau)}$ where $T$ is a one point set which maps to $e_\beta$ under $\tau$. Clearly $\rho_* = \{\rho_\epsilon\}_\epsilon$ and $\rho_*' = \{\rho_\epsilon'\}_\epsilon$ are maps of the inverse systems involved and $\rho_* \rho_*' = 1$. On the other hand, it is easily seen that $\rho_\epsilon' \rho_\epsilon([1_M]) = [1_M]$ for any geometric module $M$. Since (as remarked above) any $[\xi] \in GK_0(B;p)_\epsilon$ maps to an element of the form $[1_M]$ when control is relaxed to $3\epsilon$, it follows that $\rho_\epsilon' \rho_\epsilon$ and $1$ are coequalized by the map $GK_0(B;p)_\epsilon \to GK_0(B;p)_{3\epsilon}$. Thus $\rho_*$ is a stable isomorphism of inverse systems. Therefore, $\{GK_0(B;p)_\epsilon\}_\epsilon$ satisfies the Mittag-Leffler condition, so its $\lim^1$ vanishes. The long exact $\lim^*$-sequence induced by (**??**) proves the rest of the results claimed for $K_0$.

The second part of the proposition follows in the same way once we have established the Mittag-Leffler condition for the system $\{GK_1(B;p)_\epsilon\}$. To do so, it suffices to show that $[\nu] \in \text{Im}[GK_1(B;p)_\delta \to GK_1(B;p)_{2\epsilon}]$ for every geometric $\epsilon$-automorphism $\nu : M \to M$, every $\epsilon > 0$, and every positive $\delta < \epsilon$. By Lemma **??**, we can write $\nu$ as a composition, $\nu = \nu_{m-1} \cdots \nu_2 \nu_1$, where $\nu_i$ is a geometric $\delta$-isomorphism $M_i \to M_{i+1}, i = 1, 2, \cdots, m-1$, and $M_1 = M_m = M$. It follows easily that each sub-composition $\nu_j \nu_{j-1} \cdots \nu_i$ with $1 \le i \le j \le m-1$ is a geometric $2\epsilon$-isomorphism. Therefore, it suffices to prove the following lemma.

**Lemma 3.6.** *Let $\varphi_i : M_i \to M_{i+1}$ be a $\delta$-isomorphism for $i = 1, 2, \cdots m-1$. Assume that $\epsilon \ge \delta$ and that each composition $\varphi_{i,j} = \varphi_{j-1}\varphi_{j-2} \cdots \varphi_i : M_i \to M_j$ is a $2\epsilon$-isomorphism. Let $\varphi = \varphi_{1,m} : M_1 \to M_m = M_1$. Then the $\delta$-automorphism $\Phi$ defined below has $[\varphi] = [\Phi] \in K_1(B;p)_{2\epsilon}$. In particular, $[\varphi] \in Im[K_1(B;p)_\delta \to K_1(B;p)_{2\epsilon}]$.*

**Proof:** Let

$$(28) \qquad \Phi = \begin{pmatrix} 0 & \varphi_{m-1} \\ d & 0 \end{pmatrix} : [M_1 \oplus \cdots] \oplus M_{m-1} \to M_1 \oplus [\cdots \oplus M_{m-1}],$$

where $d$ is the diagonal morphism

$$\text{diag}(-\varphi_1, \cdots, -\varphi_{m-2}) : M_1 \oplus \cdots \oplus M_{m-2} \to M_2 \oplus \cdots \oplus M_{m-1}.$$

If each $\varphi_i$ were an honest isomorphism (in some category), then a trivial exercise would show that $\Phi$ could be reduced to the form

$$\varphi \oplus 1 : M_1 \oplus [M_2 \oplus \cdots \oplus M_{m-1}] \to M_1 \oplus [M_2 \oplus \cdots \oplus M_{m-1}],$$

by elementary column operations, i.e., by right multiplication by elementary matrices. To do the actual proof, one has to write out such a reduction; to observe that each of the elementary matrices will represent a $2\epsilon$-elementary automorphism; and to check that the homotopies which occur because we deal with $2\epsilon$-isomorphisms rather than real isomorphisms are $4\epsilon$-homotopies. We leave the details to the reader. $\square$

**3.3. Ranicki and Yamasaki's $\epsilon$-controlled K-theory.** In [**?**], two inverse systems of groups, $\{\tilde{K}_0(B, p, 0, \epsilon)\}_\epsilon$ and $\{Wh(B, p, 1, \epsilon)\}_\epsilon$, are introduced. In this section we compare these systems to those defined by us. The following proposition is essential for our treatment of the lower groups in the next section.

**Proposition 3.7.** *The inverse systems $\{\tilde{K}_0(B; p)_\epsilon\}_\epsilon$ and $\{Wh(B; p)_\epsilon\}_\epsilon$ are stably isomorphic to $\{\tilde{K}_0(B, p, 0, \epsilon)\}_\epsilon$ and $\{Wh(B, p, 1, \epsilon)\}_\epsilon$, respectively.*

**Proof:** We note that our *geometric modules*, *$\epsilon$-morphisms* and *$\epsilon$-homotopies* coincide with the *geometric modules*, the *geometric morphisms of radius $\epsilon$* and the *$\epsilon$ homotopies* in [**?**][6]. In [**?**], an *$\epsilon$-projection* $\xi : M \to M$ is required to have bound $\epsilon$ and satisfy $\xi^2 \simeq_\epsilon \xi$ while we have required an *$\epsilon$-idempotent* to have bound $\epsilon$ and satisfy $\xi^2 \simeq_{2\epsilon} \xi$. Thus, each $\epsilon$-idempotent is a $2\epsilon$-projection. As such, it can be viewed as a 0-dimensional $2\epsilon$ projective chain complex $(M, \xi)$ in the sense of the Definition following Proposition 2.2 in [**?**]. Moreover, if $\xi \sim_\epsilon \eta$ so that (**??**) and (**??**) hold, then $\psi$ and $\varphi$ show that $(M, \xi)$ and $(N, \eta)$ are 0-stable $2\epsilon$ chain equivalent (in the sense of the definitions preceeding Proposition 3.1 in [**?**]). Thus they represent the same element of $\tilde{K}_0(B, p, 0, 2\epsilon)$, and there is a well defined homomorphism $\rho_\epsilon : K_0(B; p)_\epsilon \to \tilde{K}_0(B, p, 0, 2\epsilon)$ with $\rho_\epsilon[\xi] = [M, \xi]$. The collection $\{\rho_\epsilon\}_\epsilon$ is compatible with the relaxation of control, so one has a short exact sequence of inverse systems

$$(29) \qquad 0 \to \{\mathrm{Ker}\rho_\epsilon\}_\epsilon \to \{K_0(B; p)_\epsilon\}_\epsilon \to \{\mathrm{Im}\rho_\epsilon\}_\epsilon \to 0.$$

If $\xi : M \to M$ is a geometric $\epsilon$-idempotent, then (cf. the proof of Proposition **??**) $\xi$ can be assumed $\simeq_{3\epsilon} 1_M$, so $(M, \xi)$ represents 0 in $\tilde{K}_0(B, p, 0, 3\epsilon)$. Therefore,

$$(30) \qquad \mathrm{Im}[GK_0(B; p)_\epsilon \subseteq K_0(B; p)_\epsilon \to K_0(B; p)_{\frac{3\epsilon}{2}}] \subseteq \mathrm{Ker}\rho_{\frac{3\epsilon}{2}}.$$

Conversely, if $\xi : M \to M$ has $[\xi] \in \mathrm{Ker}\rho_\epsilon$, then, by Proposition 3.1 of [**?**], $(M, \xi) \oplus (E, 1)$ is $3\epsilon$ chain equivalent to $(F, 1)$ for suitable modules $E$ and $F$. As noticed in [**?**] (immediately before Proposition 2.3), this implies that $(F, 1)$ and $(M \oplus E, \xi \oplus 1)$ are $3\epsilon$ isomorphic. But then $[\xi] = [1_F] - [1_E] \in K_0(B; p)_{3\epsilon}$, so $[\xi] \in GK_0(B; p)_{3\epsilon}$, and we have shown that

$$(31) \qquad \mathrm{Im}[\mathrm{Ker}\rho_\epsilon \subseteq K_0(B; p)_\epsilon \to K_0(B; p)_{3\epsilon}] \subseteq GK_0(B; p)_{3\epsilon}.$$

---

[6]Actually, here we are slightly inaccurate, as Ranicki and Yamasaki's geometric morphisms have $3(s, \omega, t) + 2(s, \omega, t) \neq 5(s, \omega, t)$. Since we always pass to homotopy classes, the inaccuracy is harmless.

The equations (**??**) and (**??**) imply that $\{\mathrm{Ker}\rho_\epsilon\}_\epsilon$ and $\{GK_0(B;p)_\epsilon\}_\epsilon$ are stably isomorphic subsystems of $\{K_0(B;p)_\epsilon\}_\epsilon$. It follows that the inverse systems $\{\tilde{K}_0(B;p)_\epsilon\}_\epsilon$ and $\{\mathrm{Im}\rho_\epsilon\}_\epsilon$ are also stably isomorphic, and we can finish our study of the $K_0$-case by showing that $\{\mathrm{Im}\rho_\epsilon\}_\epsilon$ is stably isomorphic to $\{\tilde{K}_0(B,p,0,\epsilon)\}_\epsilon$. A map from one to the other is given by the inclusions $\mathrm{Im}\rho_\epsilon \subseteq \tilde{K}_0(B,p,0,2\epsilon)$, and this is a stable isomorphism of inverse systems because of the inclusions $\mathrm{Im}[\tilde{K}_0(B,p,0,2\epsilon) \to \tilde{K}_0(B,p,0,4\epsilon)] \subseteq \mathrm{Im}\rho_{2\epsilon}$. The latter inclusions, in turn, hold because any $2\epsilon$ projection in the sense of [**?**] is a $2\epsilon$ idempotent in our sense.

We next turn to the $K_1$ case. Any $\epsilon$-automorphism $\varphi : M \to M$ can be viewed as a 1-dimensional, free $\epsilon$ chain complex, to be denoted $(\varphi : M \to M)$, in the sense of the Definition following Proposition 2.2 in [**?**]. As such it is strongly[7] $\epsilon$ chain contractible, so it represents an element $[\varphi : M \to M] \in Wh(X,p,1,\epsilon)$, cf. the definition preceding Proposition 4.1 in [**?**]. If $\psi : N \to N$ is another $\epsilon$-automorphism, and (**??**) holds, then the pair $(\gamma\nu, \nu)$ describes a 1-stable $2\epsilon$-simple[8] chain equivalence between $(\varphi : M \to M)$ and $(\psi : N \to N)$. Since the equivalence relation in the above mentioned definition in [**?**] is generated by 1-stable $40\epsilon$-simple chain equivalence, there results a homomorphism $\rho_\epsilon : K_1(B;p)_\epsilon \to Wh(B,p,1,\epsilon)$ given by $\rho_\epsilon([\varphi]) = [\varphi : M \to M]$, and a short exact sequence of inverse systems

$$(32) \qquad 0 \to \{\mathrm{Ker}\rho_\epsilon\}_\epsilon \to \{K_1(B;p)_\epsilon\}_\epsilon \to \{\mathrm{Im}\rho_\epsilon\}_\epsilon \to 0.$$

If $\varphi : M \to M$ is an $\epsilon$-automorphism with $\rho_\epsilon[\varphi] = 0$, then by Proposition 4.1 of [**?**], the $\epsilon$ chain complex $(\varphi : M \to M)$ is 1-stable $86\epsilon$-simple chain equivalent to the 0 complex, i.e. there are free modules $T$ and $T'$ and two $86\epsilon$ deformations, say $D_1 = f_m, f_{m-1}, \cdots, f_1$ and $D_0 = g_n, g_{n-1}, \cdots, g_1$, such that the compositions $f = f_m f_{m-1} \cdots f_1$ and $g = g_n g_{n-1} \cdots g_1$ and their inverses define $86\epsilon$ chain maps $(f,g) : (\varphi \oplus 1 : M \oplus T \to M \oplus T) \to (1 : T' \to T')$, respectively, $(f^{-1}, g^{-1}) : (1 : T' \to T') \to (\varphi \oplus 1 : M \oplus T \to M \oplus T)$. It follows that $\varphi \oplus 1_T \simeq_{172\epsilon} g_1^{-1} \cdots g_m^{-1} f_m \cdots f_1$. Since each of the factors is either a $172\epsilon$-elementary automorphism or a geometric $172\epsilon$-isomorphism, it follows from Lemma **??** that

$$(33) \qquad \mathrm{Im}[\mathrm{Ker}\rho_\epsilon \subseteq K_1(B;p)_\epsilon \to K_1(B;p)_{344\epsilon}] \subseteq GK_1(B;p)_{344\epsilon},$$

Conversely, it is easily seen that any geometric $\epsilon$-isomorphism $\xi : M \to M$ when viewed as chain complex $\xi : M \to M$ is 1-stable $\epsilon$-simple equivalent to the 0 chain complex. Therefore

$$(34) \qquad\qquad\qquad GK_1(B;p)_\epsilon \subseteq \mathrm{Ker}\rho_\epsilon$$

---

[7]Actually, "strongly" is vacuous for contractions of 1-dimensional complexes.

[8]We get $2\epsilon$ rather than $\epsilon$, because [**?**] insists that an $\epsilon$ chain map between free $\epsilon$ chain complexes commute with the boundary maps up to $\epsilon$ homotopy.

The equations (**??**) and (**??**) imply that $\{\mathrm{Ker}\rho_\epsilon\}_\epsilon$ and $\{GK_1(B;p)\}_\epsilon$ are stably isomorphic subsystems of $\{K_1(B;p)\}_\epsilon$.

As in the proof of the $K_0$ case, it now suffices to show that the system of inclusions $\{\mathrm{Im}\rho_\epsilon\}_\epsilon \subseteq \{Wh(B,p,1,\epsilon)\}_\epsilon$ is a stable isomorphism. To that end, let $[\delta : C_1 \to C_0]$ be a typical generator for $Wh(B,p,1,\epsilon)$, i.e., $\delta : C_1 \to C_0$ is a free $\epsilon$ chain complex which is $\epsilon$ chain contractible. Since an $\epsilon$ chain contraction is precisely the same as an $\epsilon$-inverse for $\delta$, Lemma **??** supplies a geometric $\epsilon$-isomorphism $\nu : C_1 \to C_0$. Since $[\delta : C_1 \to C_0] = [\delta\nu^{-1} : C_0 \to C_0] = \rho_{2\epsilon}[\delta\nu^{-1}] \in Wh(B,p,1,2\epsilon)$, it follows that

$$\mathrm{Im}[Wh(B,p,1,\epsilon) \to Wh(B,p,1,\epsilon] \subseteq \mathrm{Im}\rho_{2\epsilon},$$

and the proof is complete.

**3.4. Lower controlled K-theory.** In the Appendix of [**?**], Ranicki and Yamasaki introduce the *lower, reduced* $K$-theory groups $\tilde{K}_{1-j}(B,p,n-1,\epsilon)$ and $Wh_{2-j}(B,p,n,\epsilon)$ for all $\epsilon > 0$ and all $j,n \in \mathbb{N}$. They show that the inverse systems $\{Wh_{1-j}(B,p,n,\epsilon)\}_\epsilon$ and $\{\tilde{K}_{1-j}(B,p,n,\epsilon)\}_\epsilon$ are stably isomorphic for all $j,n \in \mathbb{N}$. Moreover, they (essentially) give the following inductive characterization of the lower groups in terms of the maps induced by the standard inclusion $i : (B;p) \to (B \times S^1; p \times id_{S^1})$.

**Proposition 3.8.** *For each $j \in \mathbb{N}$ and each $n \in \mathbb{N}$, there is a natural stable isomorphism of inverse systems*

$$\{Coker[i_* : \tilde{K}_{2-j}(B,p,n,\epsilon) \to \tilde{K}_{2-j}(B \times S^1, p \times id, n, \epsilon)]\}_\epsilon \to \{\tilde{K}_{1-j}(B,p,n,\epsilon)\}_\epsilon.$$

**Proof:** (Sketch) Of course, $\tilde{K}_{2-j}$ has to be interpreted as $Wh$ when $j = 1$. In that case, the result then follows from the stably exact sequence

$$0 \to Wh(B,p,n,\epsilon) \to Wh(B \times S^1, p \times id, n, \epsilon) \to \tilde{K}_0(B,p,n,\epsilon) \to 0,$$

which is established immediately after Corollary 7.3 in [**?**] as a formal consequence of the Mayer-Vietoris sequence (Theorem 6.2 of [**?**]) for the triad $B \times (S^1; S^1_+, S^1_-)$ and Lemma 7.2 of [**?**]. In the Appendix of [**?**] it is shown that Theorem 6.2 generalizes to lower, reduced $K$-theory, and it is easy to check a generalized version of Lemma 7.2. Therefore, the whole proof goes through for all $j \in \mathbb{N}$.                  $\square$

**Definition 3.9.** *For any $j \in \mathbb{N}$, we define $K_{-j}(B;p)_\epsilon$ to be the group $\tilde{K}_{-j}(B,p,1,\epsilon)$ of [**?**].*

Since the Ranicki-Yamasaki groups are independent of the choice of $n$ up to stable isomorphism (see [**?**], immediately before Corollary 3.5 and immediately before Proposition 4.8), and since the inclusion $i : (B;p) \to (B \times S^1, p \times id)$ is easily seen to induce a stable isomorphism on the inverse systems $\{GK_0(-;-)_\epsilon\}_\epsilon$, it is not hard to derive the following inductive characterization of the inverse systems of unreduced groups.

**Proposition 3.10.** *For each $j \in \mathbb{N}$, there is a natural stable isomorphism of inverse systems*

$$\{Coker[i_* : K_{1-j}(B; p)_\epsilon) \to K_{1-j}(B \times S^1; p \times id)_\epsilon)]\}_\epsilon \to \{K_{-j}(B; p)_\epsilon)\}_\epsilon.$$

The morphisms $i_*$ in the above propositions have natural splittings induced by the projection $(B \times S^1; p \times id) \to (B; p)$. Therefore, one gets the following corollary, where naturality refers to morphisms $(B; p) \to (B'; p')$ in the sense of (**??**). Note that the first of the three sequences involves degrees[9] $1 - j < 1$ while the other two hold in degrees $2 - j \leq 1$. Also note that the third sequence coincides with its own reduced version because of Proposition **??**.

**Corollary 3.11.** *For any $j \in \mathbb{N}$, there are natural, naturally split, exact sequences*

$$0 \to K_{1-j}(B; p)_c \xrightarrow{i_*} K_{1-j}(B \times S^1; p \times 1)_c \to K_{-j}(B; p)_c \to 0,$$

$$0 \to \tilde{K}_{2-j}(B; p)_c \xrightarrow{i_*} \tilde{K}_{2-j}(B \times S^1; p \times 1)_c \to \tilde{K}_{1-j}(B; p)_c \to 0,$$

*and*

$$0 \to \lim{}^1 K_{2-j}(B; p)_\epsilon \xrightarrow{i_*} \lim{}^1 K_{2-j}(B \times S^1; p \times 1)_\epsilon \to \lim{}^1 K_{1-j}(B; p)_\epsilon \to 0.$$

This corollary will be used in Section **??** in order to extend the basic exact sequence to lower $K$-theory.

## 4. $K$-THEORY WITH CONTINUOUS CONTROL.

### 4.1. Recollections from [**?**]. **Warning:** The functor that we call $K_*^{cc}$ in the present paper is denoted $\widetilde{K}_*^{cc}$, or sometimes just $\widetilde{K}_*$, in [**?**]. Here we use the notation $\tilde{K}_*^{cc}$ for a quotient of $K_*^{cc}$ obtained by factoring out a certain subfunctor. In particular, for us $\tilde{K}_1^{cc}$ is a just a different name for the $cc$ Whitehead group functor, $Wh^{cc}$. Further details are given below in (**??**). $\qquad\square$

Continuously controlled $K$-theory "with constant coefficients" was introduced in [**?**]. A version with "variable coefficients" was introduced by two of the present authors in [**?**] as a functor $K_*^{cc} : \mathcal{TOP}/\mathcal{CM}^* \to \mathcal{AB}_*$ Here $\mathcal{AB}_*$ is the category of $\mathbb{Z}$-graded abelian groups while $\mathcal{TOP}/\mathcal{CM}^*$ is the category of *diagrams of holink type*. We recall that an object in this category is a pair of maps $\xi = (B \xleftarrow{p} H \xrightarrow{q} X)$ with $H$ and $X$ Hausdorff, and $B$ compact, metrizable. Also, a morphism is a triple of maps making an obvious diagram commutative.

If $j > 0$, the functor $K_j^{cc}$ arises as a composition

$$\mathcal{TOP}/\mathcal{CM}^* \xrightarrow{c} \mathcal{TOP}^{cc}/\mathcal{LC} \xrightarrow{\mathcal{GM}^{cc}} \mathcal{ADDCAT} \xrightarrow{K_j^Q} \mathcal{AB},$$

---

[9]Actually, one *can* analyze the map $i_* : GK_1(B; p) \to GK_1(B \times S^1; p \times id)$ and prove exactness of the first sequence also for $1 - j = 1$, but we do not need this fact here.

where the ingredients are described in detail in [**?**]. Here we recall the major facts.

We start with the category $\mathcal{TOP}^{cc}/\mathcal{LC}$ whose objects are triples $(\overline{X}, B; p)$ with $\overline{X}$ a locally compact, Hausdorff space, $B \subseteq \overline{X}$ a compact subspace, and $p$ a map $p : E \to X = \overline{X} - B$. Even though $E$ is a very important part of such an object, it is part of the notation only insofar as any map determines its domain. Morphisms $f : (\overline{X}, B; p) \to (\overline{X}', B'; p')$ are maps $f : E \to E'$ admitting a suitable extension $\overline{f} : \overline{E} \to \overline{E}'$ where $\overline{E} = E \amalg B$ with a teardrop topology.

The *external cone* functor $c : \mathcal{TOP}/\mathcal{CM}^* \to \mathcal{TOP}^{cc}/\mathcal{LC}$ maps $\xi = (B \xleftarrow{p} H \xrightarrow{q} X)$ to the triple $(v(B), B; \mathring{v}(p))$ where $v(B)$ is the standard cone on $B$, and $\mathring{v}(p) : cyl(q) \to \mathring{v}(B) = v(B) - B$ is the map from the open mapping cylinder of $q$ (with a teardrop topology; see [**?**]) to the open cone on $B$, induced by $p$.

To each $(\overline{X}, B; p)$ we associate an additive subcategory $GM^{cc}(\overline{X}, B; p) \subseteq GM(E)$, and a corresponding additive homotopy category $\mathcal{GM}^{cc}(\overline{X}, B; p)$, by imposing the following conditions on objects, morphisms and homotopies in $GM(E)$. Each object $(S, \sigma)$ must have $p\sigma : S \to B$ proper. For each morphism $\varphi$, the family $\mathcal{F}(\varphi)$ defined in (**??**) must have $p_*\mathcal{F}(\varphi)$ *continuously controlled* at $B$ (briefly, *cc*) in the sense of [**?**]. Similarly, $p_*\mathcal{F}(\Phi)$ must be *cc* for each homotopy $\Phi$. If $\varphi$ is a morphism in $GM^{cc}(\overline{X}, B; p)$, then its image in $\mathcal{GM}^{cc}(\overline{X}, B; p)$ will be denoted $cls(\varphi)$.

For each additive category $\mathcal{A}$ and each $j \in \mathbb{N}$, we denote by $K_j^Q(\mathcal{A})$ the Quillen $K_j$-group of $\mathcal{A}$, as defined, e.g., in [**?**].

The above three functors suffice to define $K_j^{cc}$ for $j \in \mathbb{N}$, but in order to define $K_{1-j}^{cc}$ for $j \in \mathbb{N}$, we also need the suspension functor $\Sigma : \mathcal{TOP}/\mathcal{CM}^* \to \mathcal{TOP}/\mathcal{CM}^*$ (which is called $S$ in [**?**]). It has

$$\Sigma(\xi) = (\Sigma(B) \longleftarrow Cyl(q) \cup_H Cyl(q) \longrightarrow X).$$

where once again we refer to [**?**] for further details if the reader needs such.

Finally then, the complete definition as given in [**?**] can be summarized as follows.

**Definition 4.1.** *Let* $j \in \mathbb{N}$, *and* $\xi \in \mathcal{TOP}/\mathcal{CM}^*$. *If* $c$ *and* $\Sigma$ *denote the external cone and the suspension, respectively, then*

$$K_j^{cc}(\xi) = K_j^Q(\mathcal{GM}^{cc}(c(\xi))), \ \text{and} \ K_{1-j}^{cc}(\xi) = K_1^{cc}(\Sigma^j(\xi)).$$

**Remarks: 1)** We do not know whether there is a filtration on the additive category $\mathcal{GM}^{cc}(c(\xi))$ such that Pedersen and Weibels lower $K$-theoretic groups of the *filtered* additive category $\mathcal{GM}^{cc}(c(\xi))$, as defined in [**?**], coincide with the groups $K_{1-j}^{cc}(\xi)$ above.

**2)** An equivalent definition takes $K_0^{cc}(\xi) = K_0^Q(\widehat{\mathcal{GM}^{cc}}(c(\xi)))$ where the symbol $\widehat{\phantom{x}}$ indicates idempotent completion.                                                                    $\square$

From the object $\xi = (B \leftarrow H \to X)$ as above, we construct two new objects

$$\xi \times S^1 = (B \times S^1 \xleftarrow{p \times id} E \times S^1 \xrightarrow{q \circ pr} X) \quad \text{and} \quad \xi_+ = (B_+ \xleftarrow{p \amalg c_X} E \amalg X \xrightarrow{q \vee id} X),$$

where + denotes the addition of a base point (in [**?**], + was put as a superscript). Proposition 8.1 of [**?**] states that the lower $K^{cc}$-theory of $\xi_+$ can be computed inductively as follows.

**Proposition 4.2.** *For each $j \in \mathbb{Z}$ and each $\xi \in \mathcal{TOP}/\mathcal{CM}^*$, there is a natural, naturally split, short exact sequence of abelian groups*

$$0 \to K_j^{cc}(\xi_+) \xrightarrow{i_*} K_j^{cc}((\xi \times S^1)_+) \to K_{j-1}^{cc}(\xi_+) \to 0,$$

*where $i_*$ is induced by the obvious inclusion $i : \xi \to \xi \times S^1$ given by a base point in $S^1$.*

If $X$ is connected, then we shall write $K_*(X)$ for the $K$-theory (including the lower groups) of the ring $R\pi_1(X)$. More generally, $K_*(X)$ will denote the direct sum of such groups, one for each path component of $X$. By Corollary 9.4 of [**?**], these groups are naturally isomorphic to $K_*^{cc}(\emptyset_\xi)$, so the second exact sequence obtained in Section 8 of [**?**] takes the following form.

**Proposition 4.3.** *There is a natural exact sequence*

$$\cdots \to K_j(X) \to K_j^{cc}(\xi) \to K_j^{cc}(\xi_+) \to K_{j-1}(X) \to \cdots, \qquad j \in \mathbb{Z},$$

*for any $\xi \in \mathcal{TOP}/\mathcal{CM}^*$.*

**4.2. The reduced functors** $Wh^{cc}$ **and** $\tilde{K}_{2-j}^{cc}$, $j \in \mathbb{N}$. We shall also consider *reduced* groups $\tilde{K}_{2-j}^{cc}(\xi)$, $j \in \mathbb{N}$. For these, the ground ring is supposed to be $\mathbb{Z}$. Moreover, in order to identify - in Proposition **??** below - the difference between the reduced and the unreduced groups, we shall assume that each path component in each of the spaces that appear in $\xi = (B \leftarrow H \to X)$ is open. We do not know whether this extra assumption is really needed.

By definition,

$$(35) \qquad \tilde{K}_{2-j}^{cc}(\xi) = K_{2-j}^{cc}(\xi)/F_{2-j}(\xi),$$

where the subgroup $F_{2-j}(\xi)$ is given as follows. For $j = 1$, $F_1(\xi)$ is is the subgroup generated by all $[cls(\varphi)]$ where the automorphism $cls(\varphi) : (S, \sigma) \to (S, \sigma)$ in $\mathcal{GM}^{cc}(c(\xi))$ has a representative $\varphi : (S, \sigma) \to (S, \sigma)$ in $GM^{cc}(c(\xi))$ which is a geometric in the sense of (**??**). For $j = 2$, $F_0(\xi)$ is the subgroup generated by all elements of the form $[S, \sigma]$, with $(S, \sigma) \in \mathcal{GM}^{cc}(c(\xi)) \subseteq \widehat{\mathcal{GM}^{cc}}(c(\xi))$, cf. Remark **2)** following Definition **??**. Finally, for $j > 2$, $F_{2-j}(\xi) = 0$, i.e., in strictly negative degrees $K_*^{cc}$ and $\tilde{K}_*^{cc}$ coincide.

We have the following partial computation of $F_1(\xi)$ and $F_0(\xi)$. In it, for any set $S$, $\mathbb{Z}S$ is the free abelian group generated by $S$, and $[S^1, X]$ is the set of free homotopy classes of maps $S^1 \to X$.

**Proposition 4.4.** *Let $\xi = (B \xleftarrow{p} H \xrightarrow{q} X) \in \mathcal{TOP}/\mathcal{CM}^*$. With the above assumption on path components, there is a natural exact sequence*

$$Coker(q_* : \mathbb{Z}[S^1, H] \to \mathbb{Z}[S^1, X]) \xrightarrow{\alpha_1} F_1(\xi) \xrightarrow{\beta} \mathbb{Z}\pi_0(H) \xrightarrow{q_*} \mathbb{Z}\pi_0(X) \xrightarrow{\alpha_0} F_0(\xi) \longrightarrow 0.$$

An outline of a proof is given below. First, however, here are a couple of applications. By specializing to objects of the form $\xi_+$, we get the following.

**Corollary 4.5.** *For any $\xi \in \mathcal{TOP}/\mathcal{CM}^*$, $F_0(\xi_+) = 0$ for $j > 1$, and there is a natural isomorphism $F_1(\xi_+) \to \mathbb{Z}\pi_0(H)$.*

A straightforward application of Propositions **??**, **??** and **??** immediately reveal that the exact sequences of the latter two propositions have the following reduced versions.

**Proposition 4.6.** *Let $\xi \in \mathcal{TOP}/\mathcal{CM}^*$ and $j \in \mathbb{N}$. There are a natural, naturally split, short exact sequence*

$$0 \to \tilde{K}^{cc}_{2-j}(\xi_+) \to \tilde{K}^{cc}_{2-j}(\xi_+) \to \tilde{K}^{cc}_{1-j}(\xi_+) \to 0$$

*and a natural exact sequence*

$$Wh(X) \to Wh^{cc}(\xi) \to Wh^{cc}(\xi_+) \to \tilde{K}_0(X) \to \tilde{K}^{cc}_0(\xi) \to \cdots.$$

**Proof of Proposition ?? - Outline:**

We start by *defining the maps $\alpha_0, \beta$, and $\alpha_1$*. For any point $y \in cyl(q) = X \amalg H \times (0, 1)$, we denote by $M_y$ the geometric module having one generator which sits at $y$, i.e., $M_y = (\{1\}, \mu_y)$ where $\mu_y(1) = y$. If $\omega$ is path in $cyl(q)$ from $y_1$ to $y_0$, then the corresponding geometric morphism is denoted $cls(\omega) : M_{y_0} \to M_{y_1}$. With this notation, for $[x] \in \pi_0 X$, we let

$$(36) \qquad\qquad \alpha_0([x]) = [M_x] \in F_0(\xi).$$

Let $\varphi : (S, \sigma) \to (S, \sigma)$ be a geometric isomorphism representing a typical element $z = [cls(\varphi)] \in F_1(\xi)$. Thus $\varphi$ consists of a permutation $\pi : S \to S$, a collection of paths $\{\omega_s | s \in S\}$ where $\omega_s$ joins $\sigma(s)$ to $\sigma(\pi(s))$, and a collection of units $\{u_s = \pm 1 | s \in S\}$, cf. (**??**). Let $A$ be a $\pi$-orbit in $S$. To $A$ there corresponds a composable sequence of paths. Let the composite path be $\omega_A$. If $A$ is finite, this is an oriented loop $\omega_A : S^1 \to cyl(q)$. If A is infinite, it is an oriented, infinite path $\omega_A : \mathbb{R} \to cyl(q)$. In both cases, it is defined only up to an oriented reparametrization. We may decompose $S$ as $S = S_{fX} \amalg S_{iX} \amalg S_H$ where $S_{fX}$ is the union of all finite orbits $A$, for which $\omega_A$ meets $X$, $S_{iX}$ is the union of those infinite orbits $A$ for which $\omega_A$ meets $X$, and $S_H$ is the union of those (finite or infinite) orbits for which $\omega_A$ avoids $X$. The local finitenes conditions imposed on morphisms guarantee that $S_{iX}$ and $S_{fX}$ consist of finitely many orbits. For later use, we note that $z$ splits similarly

$$(37) \qquad\qquad z = z_{fX} + z_{iX} + z_H \in F_1(\xi).$$

Let $A$ be an orbit in $S_{iX}$. The oriented path $\omega_A$ must stay away from $X$ near $\pm\infty$. Therefore, it determines two path components of $H$ which we shall denote $[\omega_A(\pm\infty)]$. In terms of these, we let

$$(38) \qquad \beta(z) = \sum_A ([\omega_A(\infty)] - [\omega_A(-\infty)]) \in \mathbb{Z}\pi_0 H,$$

where the sum extends over all orbits in $S_{iX}$.

Ultimately, generators $[a]$ in the domain of $\alpha_1$ are represented by maps $a : S^1 \to X$. We let $1 \in S^1$ be the canonical base point. Then $a$ is a path from $a(1)$ to itself, so it determines a geometric automorphism $cls(a) : M_{a(1)} \to M_{a(1)}$, and we let

$$(39) \qquad \alpha_1([a]) = [cls(a)] \in F_1(\xi).$$

It is very easy to show that $\alpha_0$ *is well defined*. To see that it is *onto*, we note that any $M = (S, \sigma) \in \mathcal{GM}^{cc}(c(\xi))$ splits as $M_X \oplus M'_X$ where $M_X = (\sigma^{-1}(X), \sigma|\sigma^{-1}(X))$ and $M'_X$ is the complement. Clearly, $[M_X]$ is in the image of $\alpha_0$, so we may assume that $M = (S, \sigma) = M'_X$. Then each $\sigma(s)$ has the form $\sigma(s) = (w_s, t_s)$, with $w_s \in H$ and $0 < t_s < 1$, and we may define another object $N = (S \times \mathbb{N}, \Sigma)$ by letting $\Sigma(s, k) = (w_s, 1 - 2^{-k}(1 - t_s))$. The obvious paths from $\sigma(s)$ to $\Sigma(s, 1)$, and from $\Sigma(s, k)$ to $\Sigma(s, k+1)$, define an isomorphism $N \to N \oplus M$. This proves that $[M] = 0 \in F_0(\xi)$, and finishes the proof.

Next, we turn to *exactness at* $\mathbb{Z}\pi_0(X)$. If $w \in H$, then the geometric module $N_w = (\mathbb{N}, \nu_w)$ with $\nu_w(k) = (w, 1 - 2^{-k}) \in H \times (0, 1)$ admit an obvious geometric (and $cc$) isomorphism $M_{q(w)} \oplus N_w \to N_w$ consisting of paths in the $(0,1)$ direction. Therefore, $Im q_* \subseteq Ker \alpha_0$. To get the opposite inclusion, we need the following lemma, a proof of which can be obtained by adaptation from that of Theorem 2.6, pp. 140ff of [?].

**Lemma 4.7.** *If two geometric modules are isomorphic in* $\mathcal{GM}^{cc}(c(\xi))$*, then there is a cc geometric isomorphism between them.*

A typical element of $\mathbb{Z}[\pi_0(X)]$ is a difference $a - b$ with $a$ and $b$ of the form $a = \sum_{i=1}^n [x_i]$, $b = \sum_{j=1}^m [y_j]$ for certain points $x_i$ and $y_j$ in $X$. We have to show that $a - b \in \mathrm{Im}(q_*)$ under the assumption that $\alpha_0(a) = \alpha_0(b)$. We let $S = \{s_1, s_2, \cdots, s_n\}$, $T = \{t_1, t_2, \cdots, t_m\}$ and define $\sigma$ and $\tau$ to have $\sigma(s_i) = x_i$ and $\tau(t_j) = y_j$. The assumption implies the existence of a geometric module $(R, \rho)$ and an isomorphism in $\mathcal{GM}^{cc}(c(\xi))$ $cls(\varphi) : (S, \sigma) \oplus (R, \rho) \longrightarrow (T, \tau) \oplus (R, \rho)$. Let $R' = \rho^{-1}(X)$ and $(R, \rho) = (R', \rho') \oplus (R'', \rho'')$. Also, let $(S', \sigma') = (S, \sigma) \oplus (R', \rho')$, and define $(T', \tau')$ similarly. If $R' = \{r_1, r_2, \cdots, r_p\}$ and $\rho'(r_k) = z_k \in X$, and we set $c = \sum_{k=1}^p [z_k]$, then we may replace $a$ and $b$ with $a' = a + c$ and $b' = b + c$, respectively. In other words, we may assume that $R'$ is empty, so that $\rho$ maps into $H \times (0, 1)$. The lemma lets us assume that the morphism $\varphi : (S, \sigma) \oplus (R, \rho) \to (T, \tau) \oplus (R, \rho)$ which represents the given isomorphism is geometric, i.e., it establishes a bijection $\beta : S \amalg R \to T \amalg R$ and supplies a path between corresponding elements (or rather their images in $cyl(q)$ under $\sigma \vee \rho$ or $\tau \vee \rho$). Now we choose subsets $R_1, R_2 \subseteq R$, so that $T \cup \beta(S) = T \amalg R_2$ and $S \cup \beta^{-1}(T) = S \amalg R_1$. Then $\beta$ restricts to give a bijection $\beta' : S \amalg R_1 \to T \amalg R_2$. By restricting attention to the paths corresponding to $S \amalg R_1$ and $T \amalg R_2$, we finally arrive at a geometric isomorphism between $(S, \sigma) \oplus (R_1, \rho_1)$ and $(T, \tau) \oplus (R_2, \rho_2)$, and a moment's reflection shows that this suffices to show that $a - b \in \mathrm{Im}(q_*)$ as desired.

We go on to show that $\beta$ *is well defined*. There is a morphism

$$\eta = (\eta_B, \eta_H, \eta_X) : \xi \to \pi_0(\xi) = (\pi_0(B) \longleftarrow \pi_0(H) \longrightarrow \pi_0(X)),$$

where $\eta_?$ is the map from a space to its (by our assumption, *discrete*) quotient space of path components. By naturality, it suffices to deal with the case of $\pi_0(\xi)$, i.e., we may assume that $X$, $H$ and $B$ are discrete. Then the groups involved are direct sums of similar groups indexed on the points of $X$, so we can further assume that $X$ is just one point. Also, since $B$ is finite, the $cc$ control is unchanged, if we replace $B$ by one point. If $H$ is infinite, the groups are direct limits of similar groups corresponding to the finite subsets of $H$ so we can take $H$ to be finite. After that is done, once again we can observe that the $cc$ condition is unchanged if we replace $H \to \{*\}$ by $1 : H \to H$. In summary, we have to show that $\beta$ is well defined when $\xi = (H \xleftarrow{1} H \to \{*\})$ with $H$ finite and discrete. If $H = \emptyset$, there is nothing to prove. If $H \neq \emptyset$, Theorem 9.1 of [?] allows us to use $bc$ $K$-theory instead of $K^{cc}$, and Theorem II of [?] shows that $F_1(\xi)$ is free abelian of rank $|H| - 1$ , generated by infinite shifts, and the rest of the proof is easy.

*Exactness at* $\mathbb{Z}\pi_0 H$ is shown as follows. For each orbit $A$ in $S_{iX}$, the path $\omega_A$ meets $X$, so $q_*[\omega_A(\infty)] = q_*[\omega_A(-\infty)]$. Therefore, $\mathrm{Im}\beta \subseteq \mathrm{Ker}q_*$. To get the opposite inclusion, we note that the kernel of $q_*$ is generated by elements of the form $[w_+] - [w_-]$ where $w_\pm \in H$ and there is a path $\omega$ from $\omega(-1) = q(w_-)$ to $\omega(1) = q(w_+)$ in $X$. Let $L = (\mathbb{Z}, \lambda)$ where $\lambda(k) = (w_+, 1 - 2^{-k})$ for $k > 1$, $\lambda(k) = (w_-, 1 - 2^k)$ for $k < -1$, and $\lambda(k) = \omega(k)$ for $k = -1, 0, 1$. An obvious collection of paths in the $(0, 1)$ direction combine with $\omega$ to define a geometric automorphism of $L$ which maps to $[w_+] - [w_-]$ under $\beta$.

Clearly $[a] \mapsto [cls(a)] \in F_1(\xi)$, for $a : S^1 \to X$, defines a homomorphism $\mathbb{Z}[S^1, X] \to F_1(\xi)$, and to see that $\alpha_1$ *is well defined*, we only have to show that $[cls(q \circ b)] = 0$ for any $b : S^1 \to H$. The proof of this proceeds as follows. For $n \in \mathbb{N}$, let $\zeta_n$ be a primitive $2^{nth}$ root of unity, and define the geometric module $M_n$ to be $M_n = (\{0, 1, 2, \cdots, 2^n - 1\}, \sigma_n)$, where $\sigma_n(i) = [1 - 2^{-n}, b(\zeta_n^i)] \in cyl(q)$. The restrictions of $b$ to the arcs from $\zeta_n^i$ to $\zeta_n^{i+1}$ for $i = 0, 1, 2, \cdots 2^n - 1$, define a geometric isomorphism $\mu_n : M_n \to M_n$. Let $\mu_0 : M_0 \to M_0$ be $cls(q \circ b) : M_{q \circ b(0)} \to M_{q \circ b(0)}$. It is easily seen that the infinite direct sums $\mu_{even} = \bigoplus_{i \geq 0} \mu_{2i}$, $\mu'_{even} = \bigoplus_{i > 0} \mu_{2i}$, and $\mu_{odd} = \bigoplus_{i \geq 0} \mu_{2i+1}$, make sense. Moreover an easy subdivision and conjugation argument along the lines of the proof of Lemma 4.3 in [**?**] gives the following identity

$$[cls(\mu_{even})] = [cls(\mu_{odd})] = [cls(\mu'_{even})] \in F_1(\xi),$$

which implies that $[cls(\mu_0)] = 0 \in F_1(\xi)$, as desired.

The proof of *exactness at* $F_1(\xi)$ relies on a closer analysis of the splitting in (**??**). The necessary ideas can be dug out of [**?**], so we shall be content with these hints: In the splitting, the elements $z_{fX}$ are precisely those in the image of $\alpha_1$; this uses a version of the *subdivision* procedure (also used above) defined in Section 4 of [**?**]. Also, the elements $z_H$ vanish; this uses an Eilenberg swindle combined with subdivison and *radial push* defined as in Section 5 of [**?**]. Finally, by the definition of $\beta$ above, the elements $z_{iX}$ precisely account for the image of $\beta$.

**4.3.** *cc K-theory of compact, metrizable, tame pairs.* Let $(\overline{Y}, C)$ be a pair of spaces. We recall from [**?**] that $C$ is *tame* in $\overline{Y}$, if there is a neighborhood $N$ of $C$ in $\overline{Y}$ and a *nearly strict deformation retraction* of $\overline{Y}$ into $C$, i.e., a strict map

$$(N \times I, C \times I \cup N \times \{0\}) \to (\overline{Y}, C)$$

which is the inclusion on $N \times \{1\}$ and the projection on $C \times I$. We denote by $\mathcal{CM}^2_{tame}$ the category of *compact* pairs $(\overline{Y}, C)$ of *metrizable* spaces with $C$ *tame* in $\overline{Y}$. The morphisms are *strict* maps $f : (\overline{Y}, C) \to (\overline{Y'}, C')$, i.e. maps of pairs with $f^{-1}(C') = C$.

The *holink diagram* of Section 2 of [**?**] defines a functor $\chi : \mathcal{CM}^2_{tame} \to \mathcal{TOP}/\mathcal{CM}^*$ with

$$\chi(\overline{Y}, C) = (Y \xleftarrow{q} H = Holink(C, \overline{Y}) \xrightarrow{p} C).$$

Here $H$ is the space of strict maps $([0,1], \{0\}) \to (\overline{Y}, C)$ while $p$ and $q$ are the evaluation maps at 0 and 1, respectively.

There is also an *inclusion functor* $\iota : \mathcal{CM}^2_{tame} \to \mathcal{TOP}^{cc}/\mathcal{LC}$ with $\iota(\overline{Y}, C) = (\overline{Y}, C; id_Y)$, where $Y = \overline{Y} - C$. For a morphism $g : (\overline{Y}, C) \to (\overline{Y'}, C')$, $\iota(g) = g|Y : Y \to Y'$. Since $\overline{Y}$ is compact, the necessary conditions for $\iota(g)$ to be a morphism in $\mathcal{TOP}^{cc}/\mathcal{LC}$ as defined in [**?**], are easily verified. We also note that when we start with $(\overline{Y}, C) \in \mathcal{CM}^2_{tame}$, take $Y = \overline{Y} - C$, and give $Y \amalg C$ the teardrop topology relative to the function $id_Y \amalg id_C$, as in Section 1 of [**?**], then we do recover the space $\overline{Y}$, i.e., there is no conflict between our present use of the notation $\overline{Y}$ and that in [**?**].

**Proposition 4.8.** *Let* $j \in \mathbb{N}$ *and* $(\overline{Y}, C) \in \mathcal{CM}^2_{tame}$. *There is a natural isomorphism*

$$K^Q_j(\mathcal{GM}^{cc}(\iota(\overline{Y}, C))) \to K^{cc}_j(\chi(\overline{Y}, C)).$$

**Proof:** We shall introduce a homotopy notion in the category $\mathcal{TOP}^{cc}/\mathcal{LC}$, show that the composite functor $K^Q_j \circ \mathcal{GM}^{cc} : \mathcal{TOP}^{cc}/\mathcal{LC} \to \mathcal{AB}$ is a homotopy functor, and prove that there is a natural homotopy equivalence of functors $c \circ \chi \simeq \iota : \mathcal{CM}^2_{tame} \to \mathcal{TOP}^{cc}/\mathcal{LC}$. The desired conclusion then follows since $K^{cc}_j \circ \chi = K^Q_j \circ c \circ \chi$ for $j \in \mathbb{N}$.

The relevant homotopy notion in $\mathcal{TOP}^{cc}/\mathcal{LC}$ is obtained by declaring the *cylinder* on $(\overline{X}, B; p)$ to be

$$(\overline{X}, B; p) \times I = ([(\overline{X} \times I) \amalg B]/\sim, B; p \times id),$$

where the equivalence relation $\sim$ identifies each $\{b\} \times I$ to the point $b \in B$, and where $p \times id : H \times I \to X \times I \subseteq [(\overline{X} \times I) \amalg B]/\sim$. One then has two inclusions

$$i_\nu : (\overline{X}, B; p) \to (\overline{X}, B; p) \times I$$

($\nu = 0, 1$), and as usual, two morphisms $f_\nu : (\overline{X}, B; p) \to (\overline{Z}, D; r)$ ($\nu = 0, 1$) are called homotopic if there is a morphism $F : (\overline{X}, B; p) \times I \to (\overline{Z}, D; r)$ with $F \circ i_\nu = f_\nu$ for $\nu = 0, 1$.

In order to show that $K_j^Q \circ \mathcal{GM}^{cc}$ is a homotopy functor for each $j \in \mathbb{N}$, it suffices to establish the following lemma.

**Lemma 4.9.** *There is a natural equivalence of additive functors between additive categories*

$$i_{0*} \cong i_{1*} : \mathcal{GM}^{cc}(\overline{X}, B; p) \to \mathcal{GM}^{cc}((\overline{X}, B; p) \times I), \ (\overline{X}, B; p) \in \mathcal{TOP}^{cc}/\mathcal{LC}.$$

**Proof:** Let $\sigma : S \to H$ be an arbitrary object of $GM^{cc}(\overline{X}, B; p)$. The image of $(S, \sigma)$ under $i_{\nu*}$ is $(S, \sigma_\nu)$ with $\sigma_\nu : S \to H \times I$ having $\sigma_\nu(s) = (\sigma(s), \nu)$, ($\nu = 0, 1$, $s \in S$). The obvious paths in the $I$ direction from $\sigma_1(s)$ to $\sigma_0(s)$ define a morphism $(S, \sigma_0) \to (S, \sigma_1)$ in $GM(H \times I)$. It is not hard to check that this morphism is $cc$, and that its class in $\mathcal{GM}^{cc}((\overline{X}, B; p) \times I)$ defines the desired natural equivalence. $\square$

Finally, the desired natural homotopy equivalence, $\chi \circ c \simeq \iota$, is essentially a special case of Lemma 2.4 of [**?**]. In fact, let $P = cyl(p, q)$ be the double mapping cylinder of the diagram $\chi(\overline{Y}, C)$ with the teardrop topology described in [**?**] or in [**?**]. We apply the lemma quoted[10] with $\epsilon = diam(\overline{Y})$ to get a strict homotopy equivalence, relative to C, $\overline{\alpha} : (P, C) \to (\overline{Y}, C)$. Let $\alpha : P - C \to \overline{Y} - C$ be the homotopy equivalence obtained by restriction, and observe that $P - C$ is the domain of $\overset{\circ}{v}(p)$ in the object $\underline{c}\chi(\overline{Y}, C) = (v(C), C; \overset{\circ}{v}(p)$ while $\overline{Y} - C = Y$ is the domain of $id_Y$ in $\iota(\overline{Y}, C) = (\overline{Y}, C; id_Y)$. It is not hard to check that $\alpha$ is a morphism $c \circ \chi(\overline{Y}, C) \to \iota(\overline{Y}, C)$, and in fact a homotopy equivalence, in the category $\mathcal{TOP}^{cc}/\mathcal{LC}$. Since the naturality is clear, this finishes the proof of Proposition **?**. $\square$

Motivated by Proposition **?**, we introduce the following definition. Note that it does *not* say that $K_*^{cc}(\overline{Y}, C; id_Y)$ is obtained from the Pedersen Weibel delooping of $\mathcal{GM}^{cc}(\overline{Y}, C; id_Y)$.

---

[10]Actually, the map $P \to \overline{Y}$ given in [**?**] is usually not continuous, but this malady can be repaired by a reference to [**?**]. See also [**?**] and [**?**].

**Definition 4.10.** *If* $(\overline{Y}, C) \in \mathcal{CM}_{tame}^2$, *then we write* $K_*^{cc}(\overline{Y}, C; id_Y)$ *for the group* $K_*^{cc}(\chi(\overline{Y}, C))$. *A similar notation is used for* $Wh^{cc}$ *and* $\tilde{K}_0^{cc}$.

## 5. THE BASIC EXACT SEQUENCE.

This section is devoted to the proof of our first main theorem. In it, $\xi = (X \xleftarrow{q} E \xrightarrow{p} B)$ is an object of the category $\mathcal{TOP}/\mathcal{CM}^*$ defined in Section **??**. Also, recall that

$$\xi_+ = (X \xleftarrow{q \vee id} E \amalg X \xrightarrow{p \amalg c} B_+).$$

**Theorem 5.1.** *For each* $j \in \mathbb{N}$, *there is a natural short exact sequence*

$$0 \to \lim{}^1 K_{2-j}(B; p)_\epsilon \to K_{2-j}^{cc}(\xi_+) \to K_{1-j}(B; p)_c \to 0.$$

**Remark:** If we take $X$ to be a point, this gives the unreduced version of Theorem 1.1 of the introduction. We leave it to the reader to formulate and prove the reduced version. $\qquad\square$

**Proof:** We first note that it suffices to prove the result for $j = 1$. In fact, the inclusion $i : \{*\} \to S^1$ induces a split monomorphism from the exact sequence for $j = 1$ and $(B; p)$ to that for $j = 1$ and $(B \times S^1; p \times id)$. In view of Corollary **??** and Proposition **??**, the resulting exact sequence of cokernels gives the desired exact sequence for $j = 2$ and $(B; p)$. One proceeds inductively to get the exact sequence for each $j \in \mathbb{N}$.

Next we observe that we may assume that $X$ is a point in which case the object $\xi$ is denoted simply $p$. In fact, it is known from [**?**] that the obvious morphism $\pi : \xi \to p$ induces an isomorphism $\pi_* : K_{2-j}^{cc}(\xi_+) \to K_{2-j}^{cc}(p_+)$ for each $j \in \mathbb{N}$.

In [**?**], $K_1^{cc}(p_+)$ is defined to be $K_1$ of the additive category $\mathcal{GM}^{cc}(c(B_+), B_+; \overset{\circ}{c}(p_+))$ where $\overset{\circ}{c}(\xi_+) : \overset{\circ}{c}(E_+) \to \overset{\circ}{c}(B_+)$ is the map of open cones induced by $p_+$. In this proof it will be convenient to work with a different additive category which carries the same $K_1$ group. In particular, we consider

$$p \times id : E \times (-\infty, 0) \to B \times (-\infty, 0),$$

and note that there is a morphism in $\mathcal{TOP}^{cc}/\mathcal{LC}$ (cf. Section **??** and [**?**])

$$f : (B \times (-\infty, 0], B \times 0; p \times id) \to (c(B_+), B_+; \overset{\circ}{c}(p_+))$$

given as follows, for $(e, t) \in E \times (-\infty, 0)$

$$f(e, t) = \begin{cases} [e, t+1] \in \overset{\circ}{c}(E) \subset \overset{\circ}{c}(E_+), & \text{if } -1 < t < 0; \\ [+, 1 + t^{-1}] \in \overset{\circ}{c}(+) \subset \overset{\circ}{c}(E_+), & \text{if } -\infty < t \le -1. \end{cases}$$

**Lemma 5.2.** *The morphism* $f$ *induces a natural isomorphism*

$$f_* : K_1(\mathcal{GM}^{cc}(B \times (-\infty, 0], B \times 0; \overset{\circ}{c}(p_+))) \to K_1(\mathcal{GM}^{cc}(c(B_+), B_+; \overset{\circ}{c}(p_+))).$$

**Proof:** The commutative diagram of maps of topological spaces

$$
\begin{array}{ccc}
E \times (-\infty, -1] & \xrightarrow{\ incl\ } & E \times (-\infty, 0) \\
{\scriptstyle f}\big\downarrow & & \big\downarrow{\scriptstyle f} \\
\overset{\circ}{c}(+) & \xrightarrow{\ incl\ } & \overset{\circ}{c}(E_+)
\end{array}
$$

is a diagram of morphisms in $\mathcal{TOP}^{cc}/\mathcal{LC}$, so it induces a commutative diagram of additive functors

$$
\begin{array}{ccc}
\mathcal{GM}^{cc}(B \times (-\infty, -1], \emptyset; p \times id) & \xrightarrow{\ incl_*\ } & \mathcal{GM}^{cc}(B \times (-\infty, 0], B \times 0; p \times id) \\
{\scriptstyle f|_*}\big\downarrow & & \big\downarrow{\scriptstyle f_*} \\
\mathcal{GM}^{cc}(c(+), +; \overset{\circ}{c}(id_+)) & \xrightarrow{\ incl_*\ } & \mathcal{GM}^{cc}(c(B_+), B_+; \overset{\circ}{c}(p_+))
\end{array}
$$

An easy argument along the lines of the proof of Theorem 3.1 in [?] shows that the categories in the left hand column are flasque, so their Quillen K-theory vanishes. The horizontal functors are clearly full and faithful, and the categories in the right hand column are Karoubi filtered (see e.g., p. 355 in [?]) by those in the left hand column. Finally, since $f|E \times [\delta, 0)$ is a homeomorphism onto its image for any $\delta \in (-1, 0)$, the two Karoubi quotient categories are equivalent under the induced functor. The desired conclusion now follows by considering the homotopy fibrations from Theorem 5.3 of [?], and the induced maps connecting them. [Actually, one has to pass to the idempotent completions in order to apply the theorem quoted, but that does not change $K_1$, nor does it disturb flasque-ness]. $\qquad\square$

In the actual proof of the theorem for $j = 1$, we shall use the following abbreviations and notation.

$$
\begin{array}{llll}
GM_\epsilon & = & GM_\epsilon(B; p); & GM^{cc} & = & GM^{cc}(B \times (-\infty, 0], B \times 0; p \times id); \\
\mathcal{GM}_\epsilon & = & \mathcal{GM}_\epsilon(B; p); & \mathcal{GM}^{cc} & = & \mathcal{GM}^{cc}(B \times (-\infty, 0], B \times 0; p \times id); \\
\lim^1 K_1 & = & \lim^1 K_1(B; p)_\epsilon; & K_1^{cc} & = & K_1(\mathcal{GM}^{cc}(B \times (-\infty, 0], B \times 0; p \times id)); \\
\lim K_0 & = & \lim K_0(B; p)_\epsilon. & & &
\end{array}
$$

If $M = (S, \sigma)$, where $\sigma : S \to E \times (-\infty, 0)$, is a geometric module on $E \times (-\infty, 0)$, and $U \subseteq B \times (-\infty, 0]$, we denote by $j_U : M|U \to M$ the inclusion of the direct summand $M|U = (\sigma^{-1}(p \times id)^{-1}(U), \sigma|)$ into $M$; by $q_U : M \to M|U$ the corresponding projection; and by $\pi_U$ the idempotent $\pi_U = j_U q_U : M \to M$. Note that $j_U, q_U, \pi_U$ are considered morphisms in $GM$ rather than $\mathcal{GM}$.

The necessary connection between geometric modules on $E \times (-\infty, 0)$ and $E$ will be given by the homotopy preserving functors

(40) $$\Sigma : GM(E \times (-\infty, 0)) \to GM(E),$$

(41) $$i_{n*} : GM(E) \to GM(E \times (-\infty, 0)), \ n \in \mathbb{Z}.$$

The ("summation") functor $\Sigma$ is induced by the projection $E \times (-\infty, 0) \to E$, and $i_n : E \to E \times (-\infty, 0)$ has $i_n(e) = (e, -2^{-n})$. In general, if $\varphi : M \to N$ in $GM^{cc}$, then $\Sigma(\varphi) : \Sigma(M) \to \Sigma(N)$ will not be in $GM_\epsilon$, even for $\epsilon = \infty$ (in fact, $\Sigma(M)$ will not be finite). However, for judiciously chosen subsets $U, V \subseteq B \times (-\infty, 0)$ and a suitable $\epsilon$, one may have $\Sigma(q_V \varphi j_U) : \Sigma(M|U) \to \Sigma(N|V)$ in $GM_\epsilon$. On the other hand, $i_{n*}$ does map $GM_\epsilon$ into $GM^{cc}$, but we will often have to consider infinite direct sums $\oplus_n i_{n*}$ in which case one has to be careful that the sums make sense in $GM^{cc}$.

The composite functors satisfy

$$(42) \qquad \Sigma i_{n*} = 1_{GM(E)} \text{ and } i_{n*}\Sigma \simeq 1_{GM(E \times (-\infty, 0))},$$

where the natural homotopy consists of paths in the $(-\infty, 0)$-direction. Such paths also define natural transformations

$$(43) \qquad \nu_n : i_{n*} \to i_{(n+1)*} \text{ and } \mu_n : i_{(n+1)*} \to i_{n*},$$

and the composite transformations are homotopy equivalent to the identity transformations, i.e.,

$$(44) \qquad \nu_n \mu_n \simeq 1 \text{ and } \mu_n \nu_n \simeq 1.$$

Also in this case, when infinite direct sums of the transformations $\nu_n$ or $\mu_n$ are formed, attention has to be paid to the $cc$ condition, and the homotopy equivalence from $i_{n*}\Sigma$ to the identity is in $GM^{cc}$ only for suitable objects.

In our dealings with $\lim$ and $\lim^1$ we shall actually work over the cofinal subset $\{2^{-n} | n \in \mathbb{Z}\}$, so that an element $x \in \lim^1 K_1$ is represented by a sequence $\{[\varphi_n]\}_n \in \prod_n K_1(B; p)_{2^{-n}}$ where $\varphi_n : M_n \to M_n$ in $GM_{2^{-n}}$, $n \in \mathbb{Z}$, and there exist $\psi_n : M_n \to M_n$ in $GM_{2^{-n}}$ and homotopies $\Phi_n : \varphi_n \psi_n \simeq 1$, $\Psi_n : \psi_n \varphi_n \simeq 1$ in $GM_{2^{1-n}}$. It is easily seen that

$$(45) \qquad \varphi_\oplus = \bigoplus_{n \in \mathbb{Z}} i_{n*}(\varphi_n) : M_\oplus = \bigoplus_{n \in \mathbb{Z}} i_{n*}(M_n) \to M_\oplus$$

defines a morphism in $GM^{cc}$. Similarly, one has $\psi_\oplus$ and homotopies $\Phi_\oplus : \varphi_\oplus \psi_\oplus \simeq 1$, $\Psi_\oplus : \psi_\oplus \varphi_\oplus \simeq 1$ in $GM^{cc}$. Thus, $\varphi_\oplus$ represents an automorphism $cls(\varphi_\oplus) \in \mathcal{GM}^{cc}$, and we claim that

$$(46) \qquad i(x) = [cls(\varphi_\oplus)] \text{ defines a monomorphism } i : \lim^1 \to K_1^{cc}.$$

First note *that* $cls(\varphi_\oplus) \in \mathcal{GM}^{cc}$ is unchanged if the sequence $\varphi_n$ is replaced by a sequence $\varphi'_n$ with $\varphi_n \simeq \varphi'_n$ in $GM_{2^{-n}}$; *that* $cls(\varphi_\oplus)$ is replaced by $cls(\varphi_\oplus)cls(\gamma_\oplus)$ with $cls(\gamma_\oplus)$ elementary in $\mathcal{GM}^{cc}$, if the sequence $\varphi_n$ is replaced by a sequence $\varphi_n \gamma_n$ [$\in GM_{2^{1-n}}$] with $\gamma_n$ elementary in $GM_{2^{-n}}$; *that* a stabilization of each $\varphi_n$ (to $\varphi_n \oplus 1_{N_n}$) leads to a stabilization of $\varphi_\oplus$ (to $\varphi_\oplus \oplus 1_{N_\oplus}$); and *that* the operation preserves (direct) sums. Therefore, $\{[\varphi_n]\}_n \to [cls(\varphi_\oplus)]$ defines a homomorphism $\prod_n K_1(B; p)_{2^{-n}} \to K_1^{cc}$.

We must show that this homomorphism factors over $\lim^1 K_1$, so we suppose that $\{[\varphi_n]\}_n$ represents 0 in $\lim^1 K_1$. This means that we have an element $\{[\sigma_n]\}_n \in \prod_n K_1(B;p)_{2^{-n}}$ so that $[\varphi_n] = [\sigma_n] - [\sigma_{n+1}] \in K_1(B;p)_{2^{-n}}$, $n \in \mathbb{Z}$. Let $\sigma_n : N_n \to N_n$ have $\theta_n$ as a $2^{-n}$-inverse. Then $-[\sigma_{n+1}] = [\theta_{n+1}]$, so we may assume that $M_n = N_n \oplus N_{n+1}$, and $\varphi_n = \sigma_n \oplus \theta_{n+1} : N_n \oplus N_{n+1} \to N_n \oplus N_{n+1}$, $n \in \mathbb{Z}$. We note that $M_\oplus \neq N_\oplus \oplus N_\oplus$, but the natural transformations $\nu_n, \mu_n$ of (**??**) give rise to an isomorphism (in $\mathcal{GM}^{cc}$) $M_\oplus \to N_\oplus \oplus N_\oplus$ which conjugates $cls(\varphi_\oplus)$ into $cls(\sigma_\oplus) \oplus cls(\theta_\oplus)$. Since $cls(\sigma_\oplus)$ and $cls(\theta_\oplus)$ are inverses in $\mathcal{GM}^{cc}$, it follows that $[cls(\varphi_\oplus)] = 0 \in K_1^{cc}$, so we do have the desired factorization.

To show that $i$ is monic, we need the following lemma, the straightforward proof of which will be left to the reader.

**Lemma 5.3.** *Let $\mathcal{F}$ be a cc family of subsets of $B \times (-\infty, 0]$ and let $U_0 = B \times (-\infty, 0]$. There is a strictly increasing function $\lambda : \mathbb{N} \to \mathbb{Z}$ such that any $F \in \mathcal{F}$ which meets $U_n = B \times [-2^{-\lambda(n)}, 0]$ is completely contained in $U_{n-1}$, and has diameter $diam(F) \leq 2^{-n}$, $n \in \mathbb{N}$.*

If this lemma is applied to $\mathcal{F}(\varphi)$ for some morphism $\varphi : M \to N$, then for any $n < m \in \mathbb{N} \cup \infty$, $\varphi$ maps $M|(U_n - U_m)$ into $N|(U_{n-1} - U_{m+1})$, and the resulting morphism $\Sigma(M|(U_n - U_m)) \to \Sigma(N|(U_{n-1} - U_{m+1}))$, has bound $2^{-n}$. A similar remark holds for homotopies.

Returning to the proof of the injectivity of $i$, we assume that $x = \{[\varphi_n]\}_n$ has $[cls(\varphi_\oplus)] = 0 \in K_1^{cc}$. After a stabilization (which can easily be absorbed into the given morphisms $\varphi_n$), then in $GM^{cc}$ we have thirteen[11] elementary automorphisms $\gamma^i : M_\oplus \to M_\oplus$, $i = 1, 2, \cdots, 13$ and a homotopy $\Gamma : \varphi_\oplus \simeq \gamma = \gamma^1 \gamma^2 \cdots \gamma^{13}$. We apply Lemma **??** to $\mathcal{F} = \mathcal{F}(\Gamma) \cup \bigcup_i \mathcal{F}(\gamma^i)$ to determine $\lambda : \mathbb{N} \to \mathbb{Z}$ and the neighborhoods $U_n$ of $B$, $n \in \mathbb{N}$.

For each pair $(n_1, n_2) \in (\mathbb{N} \cup \{\infty\})^2$ with $n_1 < n_2$, let $\mathbf{n} = [n_1, n_2)$ be the corresponding half open interval. Let $M_\mathbf{n} = M_\oplus|(U_{n_1} - U_{n_2})$ with *complement* $M_\mathbf{n}^\perp$ ($= M_\oplus|(U_0 - U_{n_1}) \oplus M_\oplus|U_{n_2})$; *inclusion* $j_\mathbf{n} : M_\mathbf{n} \to M_\oplus$; *projection* $q_\mathbf{n} : M_\oplus \to M_\mathbf{n}$; and *idempotent* $\pi_\mathbf{n} = j_\mathbf{n} q_\mathbf{n} : M_\oplus \to M_\oplus$. Also let $\varphi_\mathbf{n} = q_\mathbf{n} \varphi_\oplus j_\mathbf{n} : M_\mathbf{n} \to M_\mathbf{n}$ (so $\Sigma \varphi_\mathbf{n} = \oplus_k \varphi_k$ with $\lambda(n_1) \leq k < \lambda(n_2)$). Finally, for the above $\mathbf{n} = [n_1, n_2)$, let

(47)     $$\mathbf{n}'' = [n_1, n_2 + 1) \; ; \; \mathbf{n}' = [n_1 - 1, n_2), \text{ if } n_1 > 1; \text{ and}$$

(48)     $$\mathbf{n} + k\mathbf{e} = [n_1 - k, n_2 + k), \text{ if } k \in \mathbb{Z} \text{ and } 0 < n_1 - k < n_2 + k.$$

For any $\mathbf{n}$ there is the elementary automorphism (of bound $2^{-n_1}$) $\gamma_\mathbf{n}^i = 1 + (\gamma^i - 1)\pi_\mathbf{n} : M_\oplus \to M_\oplus$. It maps $M_{\mathbf{n}+\mathbf{e}}$ into itself, is the identity on $M_{\mathbf{n}+\mathbf{e}}^\perp$, and agrees with the restriction of $\gamma^i$ on the submodule $M_{\mathbf{n}} \subseteq M_{\oplus}$. We think of $\gamma_{\mathbf{n}}^i$ as a "concentration of $\gamma^i$ on $M_{\mathbf{n}}$" although this is not literally true. If $n_2 > n_1 + 26$, then the composition

$$(49) \qquad \gamma_{\mathbf{n}} = \gamma_{\mathbf{n}-\mathbf{e}}^1 \gamma_{\mathbf{n}-2\mathbf{e}}^2 \cdots \gamma_{\mathbf{n}-13\mathbf{e}}^{13} : M_{\oplus} \to M_{\oplus}$$

has bound $2^{-n_1-13}(1+2+4+\cdots+2^{12}) < 2^{-n_1}$. Moreover, each factor maps $M_{\mathbf{n}}$ into itself and is the identity on $M_{\mathbf{n}}^{\perp}$. Hence $q_{\mathbf{n}}\gamma_{\mathbf{n}}j_{\mathbf{n}} : M_{\mathbf{n}} \to M_{\mathbf{n}}$ is a product of the 13 elementary automorphisms $q_{\mathbf{n}}\gamma_{\mathbf{n}-i\mathbf{e}}^i j_{\mathbf{n}} : M_{\mathbf{n}} \to M_{\mathbf{n}}$, so

$$(50) \qquad [\Sigma(q_{\mathbf{n}}\gamma_{\mathbf{n}}j_{\mathbf{n}})] = 0 \in K_1(B;p)_{2^{-n_1}}, \text{ if } n_2 > n_1 + 26.$$

---

[11]In [**?**] (see pp. 145-46), two of the present authors observed that 16 elementary automorphisms will do. The reduction to 13 factors can be found in [**?**], or by contemplating the proof of Proposition **??**.

We now assume that $n_2 - n_1 \geq 39$, and write $M_{\mathbf{n}} = M_{\mathbf{l}} \oplus M_{\mathbf{m}} \oplus M_{\mathbf{r}}$ where $\mathbf{l} = [n_1, n_1+13]$, $\mathbf{m} = \mathbf{n} - 13\mathbf{e}$, and $\mathbf{r} = [n_2 - 13, n_2]$. The indexing is supposed to convey the **l**eft hand, the **m**iddle, and the **r**ight hand summand of $M_{\mathbf{n}}$, respectively. We note that $M_{\mathbf{r}} = 0$ when $n_2 = \infty$. Since the length of the middle interval $\mathbf{m}$ is at least 13, $\gamma_{\mathbf{n}}$ maps the left summand $M_{\mathbf{l}}$ into $M_{\mathbf{l}} \oplus M_{\mathbf{m}}$ and the right summand $M_{\mathbf{r}}$ into $M_{\mathbf{m}} \oplus M_{\mathbf{r}}$. Also, $\gamma_{\mathbf{n}}$ coincides with $\gamma$ on $M_{\mathbf{m}}$, so

$$q_{\mathbf{n}}\gamma_{\mathbf{n}}j_{\mathbf{m}} = q_{\mathbf{n}}\gamma j_{\mathbf{m}} \simeq q_{\mathbf{n}}\varphi_{\oplus}j_{\mathbf{m}} : M_{\mathbf{m}} \to M_{\mathbf{n}}$$

(the homotopy is the reverse of $q_{\mathbf{n}}\Gamma j_{\mathbf{m}}$ and has bound $2^{-n_1}$). In terms of matrices w.r.t. the decomposition $M_{\mathbf{n}} = M_{\mathbf{l}} \oplus M_{\mathbf{m}} \oplus M_{\mathbf{r}}$ this gives us the equation

$$(51) \qquad q_{\mathbf{n}}\gamma_{\mathbf{n}}j_{\mathbf{n}} = \begin{pmatrix} \gamma_{\mathbf{n}}^{(l)} & q_{\mathbf{l}}\gamma_{\mathbf{n}}j_{\mathbf{m}} & 0 \\ * & q_{\mathbf{m}}\gamma_{\mathbf{n}}j_{\mathbf{m}} & * \\ 0 & q_{\mathbf{r}}\gamma_{\mathbf{n}}j_{\mathbf{m}} & \gamma_{\mathbf{n}}^{(r)} \end{pmatrix} \simeq \begin{pmatrix} \gamma_{\mathbf{n}}^{(l)} & 0 & 0 \\ * & \varphi_{\mathbf{m}} & * \\ 0 & 0 & \gamma_{\mathbf{n}}^{(r)} \end{pmatrix}$$

in $GM(E \times (-\infty, 0); p \times id)_{2^{-n_1}}$ (note that the center columns represent $q_{\mathbf{n}}\gamma j_{\mathbf{m}}$ and $q_{\mathbf{n}}\varphi_{\oplus}j_{\mathbf{m}}$, respectively, and that the equation *defines* $\gamma_{\mathbf{n}}^{(l)}$ and $\gamma_{\mathbf{n}}^{(r)}$).

Since the left hand side and $\varphi_{\mathbf{m}}$ are $2^{-n_1}$-isomorphisms, it easily follows that $\gamma_{\mathbf{n}}^{(l)}$ and $\gamma_{\mathbf{n}}^{(r)}$ are $2^{2-n_1}$-isomorphisms. Clearly, (**??**) and (**??**) imply that[12]

$$(52) \qquad [\Sigma\gamma_{\mathbf{n}}^{(l)}] + [\Sigma\varphi_{\mathbf{m}}] + [\Sigma\gamma_{\mathbf{n}}^{(r)}] = 0 \in K_1(B;p)_{2^{3-n_1}}, \text{ if } 39 \leq n_2 - n_1 < \infty.$$

It is easily seen that the upper left hand corner of the above matrix is unchanged if one augments $n_2$, keeping $n_1$ fixed, i.e., $\gamma_{\mathbf{n}}^{(l)} = \gamma_{\mathbf{n}''}^{(l)}$ (and, similarly, $\gamma_{\mathbf{n}}^{(r)} = \gamma_{\mathbf{n}'}^{(r)}$). Hence the definition $b_{n_1-3} = [\Sigma(\gamma_{\mathbf{n}}^{(l)})] \in K_1(B;p)_{2^{3-n_1}}$ gives an element which is independent of the choice of $n_2 \in [n_1+39, \infty)$. For any $\mathbf{n} \in \mathbb{N}$, we can apply (**??**) to

---

[12]The loss of a factor 2 in control happens because one must get rid of the off diagonal terms by an elementary operation.

$\mathbf{n} = [n+4, n+43)$ and $\mathbf{n}' = [n+3, n+43)$ to get

$$b_{n+1} + \sum_{k=\lambda(n+17)}^{\lambda(n+30)-1} [\varphi_k] + [\Sigma(\gamma_{\mathbf{n}}^{(r)})] \;\; = \;\; 0 \in K_1(B;p)_{2^{-n-1}}, \text{ and}$$

$$b_n + \sum_{k=\lambda(n+16)}^{\lambda(n+30)-1} [\varphi_k] + [\Sigma(\gamma_{\mathbf{n}'}^{(r)})] \;\; = \;\; 0 \in K_1(B;p)_{2^{-n}}.$$

If we let $c_n = -b_n + \sum_{k=n}^{\lambda(n+16)-1} [\varphi_k] \in K_1(B;p)_{2^{-n}}$, then the above formulae immediately imply that $c_n - r_{n,n+1}(c_{n+1}) = [\varphi_n] \in K_1(B;p)_{2^{-n}}$ for all $n \in \mathbb{N}$. This shows that $x = \{[\varphi_k]\}_k = 0 \in \lim^1 K_1$ and finishes the proof of the injectivity of $i$.

We go on to define the map $\partial : K_1^{cc} \to \lim K_0$. Let $y = [cls\,\varphi] \in K_1^{cc}$. Thus $\varphi : M \to M$ in $GM^{cc}$ and there exist $\psi : M \to M$ and homotopies $\Phi : \varphi\psi \simeq 1$ and $\Psi : \psi\varphi \simeq 1$ in $GM^{cc}$. We use Lemma **??** with $\mathcal{F} = \mathcal{F}(\varphi) \cup \mathcal{F}(\psi) \cup \mathcal{F}(\Phi) \cup \mathcal{F}(\Psi)$ and continue to denote the resulting neighborhoods by $U_n$. Also, the notation $j_{\mathbf{n}}, q_{\mathbf{n}}, \pi_{\mathbf{n}}, M_{\mathbf{n}}$, introduced in relation to $M_{\oplus}$ immediately before (**??**), will now be used for $M$. Thus

$$\pi_{\mathbf{n}} \varphi j_{\mathbf{k}} = \varphi j_{\mathbf{k}} \text{ whenever } \mathbf{n} \supseteq \mathbf{k} + \mathbf{e},$$

and analogous identities hold for $\psi, \Phi$, and $\Psi$. For any $\mathbf{n} = [n_1, n_2)$ with $n_2 - n_1 \geq 3$, we let

$$(53) \qquad p_{\mathbf{n}} = p_{\mathbf{n}}(\varphi, \psi) = q_{\mathbf{n}} \varphi \pi_{\mathbf{n-e}} \psi j_{\mathbf{n}} : M_{\mathbf{n}} \to M_{\mathbf{n}}$$

Since $j_{\mathbf{n}} q_{\mathbf{n}} \varphi \pi_{\mathbf{n-e}} = \varphi \pi_{\mathbf{n-e}}$, the homotopy $q_{\mathbf{n}} \varphi \pi_{\mathbf{n-e}} \Psi \pi_{\mathbf{n-e}} \psi j_{\mathbf{n}}$ shows that $p_{\mathbf{n}}$ is a homotopy idempotent of bound $2^{-n_1}$ on the module $M_{\mathbf{n}}$.

If we assume that $n_2 - n_1 \geq 6$, then we get a decomposition into a left hand, a middle, and a right hand summand as before, $M_{\mathbf{n}} = M_{\mathbf{l}} \oplus M_{\mathbf{m}} \oplus M_{\mathbf{r}}$, but this time $\mathbf{l} = [n_1, n_1+2)$, $\mathbf{m} = \mathbf{n} - 2\mathbf{e}$, $\mathbf{r} = [n_2 - 2, n_2)$. As before, the right hand summand vanishes if $n_2 = \infty$. Since $\mathbf{m}$ has length at least two, $p_{\mathbf{n}}$ maps $M_{\mathbf{l}}$ into $M_{\mathbf{l}} \oplus M_{\mathbf{m}}$ and $M_{\mathbf{r}}$ into $M_{\mathbf{m}} \oplus M_{\mathbf{r}}$. Also,

$$p_{\mathbf{n}} | M_{\mathbf{m}} = q_{\mathbf{n}} \varphi \pi_{\mathbf{n-e}} \psi j_{\mathbf{m}} = q_{\mathbf{n}} \varphi \psi j_{\mathbf{m}} \simeq q_{\mathbf{n}} j_{\mathbf{m}} = incl : M_{\mathbf{m}} \to M_{\mathbf{n}}$$

by the homotopy $q_{\mathbf{n}} \Phi j_{\mathbf{m}}$ of bound $2^{-m_1} = 2^{-(n_1+2)}$. In matrix language relative to the decomposition $M_{\mathbf{n}} = M_{\mathbf{l}} \oplus M_{\mathbf{m}} \oplus M_{\mathbf{r}}$, this means that we have a homotopy of bound $2^{-m_1}$

$$(54) \qquad p_{\mathbf{n}} = \begin{pmatrix} p_{\mathbf{n}}^{(l)} & q_{\mathbf{l}} \varphi \psi j_{\mathbf{m}} & 0 \\ * & q_{\mathbf{m}} \varphi \psi j_{\mathbf{m}} & * \\ 0 & q_{\mathbf{r}} \varphi \psi j_{\mathbf{m}} & p_{\mathbf{n}}^{(r)} \end{pmatrix} \simeq \begin{pmatrix} p_{\mathbf{n}}^{(l)} & 0 & 0 \\ * & 1 & * \\ 0 & 0 & p_{\mathbf{n}}^{(r)} \end{pmatrix}$$

in $GM(E \times (-\infty, 0); p \times id)$. It follows easily that the right hand matrix is a homotopy idempotent of bound $2^{1-n_1}$. Hence $p_{\mathbf{n}}^{(l)}$ and $p_{\mathbf{n}}^{(r)}$ (which are *defined* by the equation) are homotopy idempotents of bound $2^{1-n_1}$. Moreover, the left hand idempotent $p_{\mathbf{n}}^{(l)}$ is unchanged if one augments $n_2$, keeping $n_1$ fixed (and dually for the right hand idempotent $p_{\mathbf{n}}^{(r)}$), i.e., with notation from (??),

$$(55) \qquad p_{\mathbf{n}'}^{(r)} = p_{\mathbf{n}}^{(r)} \text{ and } p_{\mathbf{n}''}^{(l)} = p_{\mathbf{n}}^{(l)}.$$

We define $\partial : K_1^{cc} \to \lim K_0(B;p)_{2^{-n}}$ by

$$(56) \qquad \partial[cls\,\varphi] = \{\partial_n(\varphi,\psi)\}_n,$$

where

$$(57) \qquad \partial_{n-1}(\varphi,\psi) = [\Sigma p_{[n,n+k]}^{(l)}(\varphi,\psi)] - [\Sigma p_{[n,n+k]}^{(l)}(1_M,1_M)], 6 \leq k \leq \infty.$$

To check that this actually gives an element of $\lim K_0$ which is independent of the various choices we need the following identity in $K_0(B;p)_{2^{1-n_1}}$

$$(58) \qquad [\Sigma p_{\mathbf{n}}(1_M,1_M)] = [\Sigma p_{\mathbf{n}}(\varphi,\psi)] = [\Sigma p_{\mathbf{n}}^{(l)}(\varphi,\psi)] + [1_{\Sigma M \mathbf{m}}] + [\Sigma p_{\mathbf{n}}^{(r)}(\varphi,\psi)].$$

The second part of (??) follows from (??). To prove the first part, let $f = q_{\mathbf{n}}\varphi\pi_{\mathbf{n-e}}j_{\mathbf{n}} : M_{\mathbf{n}} \to M_{\mathbf{n}}$ and $g = q_{\mathbf{n}}\pi_{\mathbf{n-e}}\psi j_{\mathbf{n}} : M_{\mathbf{n}} \to M_{\mathbf{n}}$. Clearly $fg = p_{\mathbf{n}}(\varphi,\psi)$, and, since $\pi_{\mathbf{n}}\varphi\pi_{\mathbf{n-e}} = \pi_{\mathbf{n-e}}$, $gf \simeq p_{\mathbf{n}}(1_M,1_M)$ with bound $2^{-1-n_1}$.

The identity shows that $\partial_{n-1}(\varphi,\psi)$ is also given by

$$(59) \qquad \partial_{n-1}(\varphi,\psi) = -[\Sigma p_{[n,n+k]}^{(r)}(\varphi,\psi)] + [\Sigma p_{[n,n+k]}^{(r)}(1_M,1_M)], 6 \leq k < \infty.$$

It follows from (??) and (??) that $\{\partial_n(\varphi,\psi)\}_n \in \lim K_0(B;p)_{2^{-n}}$. It also follows that the choice of the neighborhoods $U_n$ satisfying Lemma ?? is immaterial. Indeed, if $\{U_n'\}_n$ is a different choice, then $\{U_n \cap U_n'\}_n$ is also a possible choice, so in a comparison we may assume that $U_n' \subseteq U_n$. But then for any $n_0$, the system $\{V_n^{(n_0)}\}_n$ which coincides with $U_n$ for $n < n_0$ and with $U_n'$ for $n \geq n_0$ is allowed and one only has to use (??) and (??) to the system $\{V_n^{(n_0)}\}_n$ for a suitable choice of $n_0$ and $k$ to finish the proof. Similarly, a modification of $\varphi$ (or $\psi$) by a homotopy or by right multiplication by an elementary $\gamma$ can be decomposed into two consecutive modifications where one does not change anything above some arbitrarily chosen $U_{n_0}$ while the other one leaves everything intact above the complement of $U_{n_0+k}$ for suitable $k$. Therefore, the above argument easily adapts to show that, indeed, (??) and (??) give a well defined homomorphism.

To show that $\partial$ is epic, let $z = \{[\alpha_n] - [\beta_n]\}_n \in \lim K_0(B;p)_{2^{-n}}$. This means that $\alpha_n : A_n \to A_n$ and $\beta_n : B_n \to B_n$ are homotopy idempotents of bound $2^{-n}$, and in $GM_{2^{-n}}$ there exist geometric modules $S_n$ and morphisms $\psi_n : A_n \oplus B_{n+1} \oplus S_n \to$ $A_{n+1} \oplus B_n \oplus S_n$ and $\varphi_n \colon A_{n+1} \oplus B_n \oplus S_n \to A_n \oplus B_{n+1} \oplus S_n$ such that we have the following homotopies of bound $2^{1-n}$

$$(60) \qquad \alpha_n^2 \quad \simeq \quad \alpha_n, \qquad\qquad \text{and } \beta_n^2 \simeq \beta_n;$$

$$(61) \qquad \varphi_n \psi_n \quad \simeq \quad \alpha_n \oplus \beta_{n+1} \oplus 1_{S_n}, \text{ and } \psi_n \varphi_n \simeq \alpha_{n+1} \oplus \beta_n \oplus 1_{S_n};$$

$$(62) \qquad \varphi_n \quad \simeq \quad \varphi_n(\alpha_{n+1} \oplus \beta_n \oplus 1_{S_n}) \simeq (\alpha_n \oplus \beta_{n+1} \oplus 1_{S_n})\varphi_n;$$

$$(63) \qquad \psi_n \quad \simeq \quad \psi_n(\alpha_n \oplus \beta_{n+1} \oplus 1_{S_n}) \simeq (\alpha_{n+1} \oplus \beta_n \oplus 1_{S_n})\psi_n.$$

We now set $T = \oplus_n i_{n*}(A_n \oplus B_n \oplus S_n)$, $T' = \oplus_n i_{n*}(A_{n+1} \oplus B_n \oplus S_n)$, and $T'' = \oplus_n i_{n*}(A_n \oplus B_{n+1} \oplus S_n)$, where $n$ ranges over $\mathbb{Z}$. These three geometric modules in $GM(E \times (-\infty, 0))$ are connected by natural transformations (cf. (??)) $\nu' \colon T' \to T$, $\mu' \colon T \to T'$, $\nu'' \colon T'' \to T$, $\mu'' \colon T \to T''$ which are pairwise inverses up to natural, cc homotopies. Put $\Phi = \oplus_n i_{n*}(\varphi_n) \colon T' \to T''$, $\Psi = \oplus_n i_{n*}(\psi_n) \colon T'' \to T'$, and $\Gamma = \oplus_n i_{n*}((1 - \alpha_n) \oplus (1 - \beta_n) \oplus 0) \colon T \to T$. Elementary computations using the homotopies from (??) – (??) show that

$$\tilde{\Phi}\Gamma \simeq 0, \Gamma\tilde{\Phi} \simeq 0, \tilde{\Psi}\Gamma \simeq 0, \Gamma\tilde{\Psi} \simeq 0, \Gamma^2 \simeq \Gamma, \tilde{\Psi}\tilde{\Phi} \simeq 1 - \Gamma, \text{ and } \tilde{\Phi}\tilde{\Psi} \simeq 1 - \Gamma,$$

where $\tilde{\Psi} = \nu'\Psi\mu'' \colon T \to T$, and $\tilde{\Phi} = \nu''\Phi\mu' \colon T \to T$. Moreover, all of these maps and homotopies are in $GM^{cc}$, so $\psi = \tilde{\Psi} + \Gamma \colon T \to T$ and $\varphi = \tilde{\Phi} + \Gamma \colon T \to T$ are cc homotopy inverses. We claim that $\partial[cls(\varphi)] = z$.

To compute $\partial_n(\varphi, \psi)$, we chose $U_n = B \times [-2^{-n}, 0]$ in Lemma ?? and $k = 6$ in (??). The homotopy idempotent $\Sigma p^{(l)}_{[n,n+6)}(\varphi, \psi)$ then lives on the module $A_n \oplus B_n \oplus S_n \oplus A_{n+1} \oplus B_{n+1} \oplus S_{n+1}$, and a simple computation gives us the following matrix form up to $2^{1-n}$-homotopies coming from (??) – (??). Here a superscript like $^{AB}$ denotes the component of a morphism which maps the relevant $B$-summand to the relevant $A$-summand.

$$\begin{pmatrix} \varphi_n^{AA}\psi_n^{AA} & 0 & \varphi_n^{AA}\psi_n^{AS} & 0 & \varphi_n^{AA}\psi_n^{AB} & 0 \\ 0 & 0 & 0 & 0 & 0 & 0 \\ \varphi_n^{SA}\psi_n^{AA} & 0 & \varphi_n^{SA}\psi_n^{AS} & 0 & \varphi_n^{SA}\psi_n^{AB} & 0 \\ 0 & 0 & 0 & 1_{A_{n+1}} & 0 & 0 \\ \varphi_n^{BA}\psi_n^{AA} & 0 & \varphi_n^{BA}\psi_n^{AS} & 0 & \varphi_n^{BA}\psi_n^{AB} + (1 - \beta_{n+1}) & 0 \\ 0 & 0 & 0 & 0 & 0 & 1_{S_{n+1}} \end{pmatrix}$$

This is the sum of two homotopy orthogonal homotopy idempotents, namely the diagonal matrix $\text{diag}(0, 0, 0, 1_{A_{n+1}}, 1 - \beta_{n+1}, 1_{S_{n+1}})$ and the product $X^t Y$ where $X = (\varphi_n^{AA}, 0, \varphi_n^{SA}, 0, \varphi_n^{BA}, 0)$, $Y = (\psi_n^{AA}, 0, \psi_n^{AS}, 0, \psi_n^{AB}, 0)$ and $^t$ denotes a transposition. These terms represent $[1_{A_{n+1}}] + [1 - \beta_{n+1}] + [1_{S_{n+1}}]$, respectively, $[X^t Y] = [Y X^t] = [(\psi_n \varphi_n)^{AA}] = [\alpha_{n+1}]$. Since $[p^{(l)}_{[n,n+6)}(1_T, 1_T)] = [1_{A_{n+1} \oplus B_{n+1} \oplus S_{n+1}}]$, it follows that $\partial_{n-1}(\varphi, \psi) = [\alpha_{n+1}] - [\beta_{n+1}] \in \acute{K}_0(B; p)_{2^{1-n}}$, so $\partial[cls\varphi] = z$ as claimed.

We leave to the reader the easy proof that $\partial i = 0$ and go on to prove that $\text{Ker}(\partial) \subseteq \text{Im}(i)$. We keep the notation from the last few lines preceding (??) and assume that $y = [cls\varphi] \in \mathrm{Ker}(\partial)$, i.e., for each $n \in \mathbb{N}$, the element $\partial_{n-1}(\varphi, \psi)$ defined in (**??**) vanishes in $K_0(B; p)_{2^{1-n}}$. The basic lemma is the following one.

**Lemma 5.4.** *There exists a universal constant $c \in \mathbb{N}$ for which the following holds. If $\partial_{n-1}(\varphi, \psi) = 0 \in K_0(B; p)_{2^{1-n}}$, then there exist a (finite) geometric module $L_{n+1}$ in $GM(B; p)$, seven $2^{c-n}$-elementary automorphisms $\gamma^i : M' = M \oplus i_{n+1*}L_{n+1} \to M \oplus i_{n+1*}L_{n+1}$, and a morphism $\varphi' \simeq \gamma^7\gamma^6\cdots\gamma^1\varphi : M' \to M'$ such that $\varphi'$ is diagonal w.r.t. the decomposition $M' = M'_{(-\infty, n+1)} \oplus M'_{[n+1, \infty)}$. Moreover, each $\gamma^i$ is the identity on $(M'_{[n, n+3)})^{\perp}$. Finally, the homotopy is bounded by $2^{c-n}$ measured in $B$ and it is stationary on $(M'_{[n, n+3)})^{\perp}$.*

**Proof:** Throughout the proof, the constant $c$ will be increased a few times without our changing its name. We take $k = \infty$ in the definition (**??**) and write $\mathbf{n} = [n, \infty)$. Then the right hand summand in the decomposition used in (**??**) vanishes, so $M = M_{\mathbf{n}}^{\perp} \oplus M_{\mathbf{l}} \oplus M_{\mathbf{m}}$ where $\mathbf{l} = [n, n+2)$ and $\mathbf{m} = [n+2, \infty)$. Up to $2^{c-n}$-homotopy, one has the following matrix presentation relative to this decomposition (compare (**??**) and (**??**)).

$$(64) \qquad \varphi\pi_{\mathbf{n-e}}\psi = \begin{pmatrix} 0 & 0 & 0 \\ 0 & p_{\mathbf{n}}^{(l)}(\varphi, \psi) & 0 \\ 0 & x & 1 \end{pmatrix}$$

Since $\varphi\pi_{\mathbf{n-e}}\psi$ is a $2^{c-n}$-idempotent, $xp_{\mathbf{n}}^{(l)}(\varphi, \psi) \simeq 0$ (with bound $2^{c+1-n}$). We let $\gamma^1 = e_{32}(x)$ where the classical index-notation for the elementary automorphism refers to the decomposition of $M$ used in (**??**). The effect of replacing $(\varphi, \psi)$ by $(\gamma^1\varphi, \psi(\gamma^1)^{-1})$ then is to replace $x$ by $0$ in (**??**) (up to $2^{1+c-n}$-homotopy). Thus we shall proceed with $x = 0$ in (**??**). Note that this first elementary $\gamma^1$ is supported on $M_{[n, n+3)}$ as claimed. The remaining six elementary automorphisms will actually be supported on $M'_{\mathbf{l}}$ (recall that $\mathbf{l} = [n, n+2)$).

Since $p_{\mathbf{n}}^{(l)}(1_M, 1_M) : M_{\mathbf{l}} \to M_{\mathbf{l}}$ is the canonical projection onto $M_{[n+1, n+2)}$, Proposition **??** gives the following translation of the vanishing of $\partial_{n-1}(\varphi, \psi)$. For some constant $c \in \mathbb{N}$ we have a geometric module $L'_n \in GM_0$; two $2^{c-n}$-morphisms

$$(65) \qquad \alpha_n : \Sigma M_{\mathbf{l}} \oplus L'_n \to \Sigma M_{[n+1, n+2)} \oplus L'_n \text{ and}$$

$$(66) \qquad \beta_n : \Sigma M_{[n+1, n+2)} \oplus L'_n \to \Sigma M_{\mathbf{l}} \oplus L'_n;$$

and four $2^{c+1-n}$-homotopies (note that two of the conditions in (**??**) are vacuous, because the second idempotent involved here is the identity)

$$(67) \qquad \beta_n\alpha_n \simeq \Sigma p_{\mathbf{n}}^{(l)}(\varphi, \psi) \oplus 1_{L'_n} \text{ and } \alpha_n\beta_n \simeq 1_{\Sigma M_{[n+1, n+2)} \oplus L'_n}; \text{ and}$$

$$(68) \qquad \alpha_n(\Sigma p_{\mathbf{n}}^{(l)}(\varphi, \psi) \oplus 1_{L'_n}) \simeq \alpha_n \text{ and } \beta_n \simeq (\Sigma p_{\mathbf{n}}^{(l)}(\varphi, \psi) \oplus 1_{L'_n})\beta_n.$$

If we stabilize $\varphi$ and $\psi$ by adding $1_{i_{n+1*}L'_n}$, the summands $L'_n$ will appear in the original data, and no stabilization will be needed. Hence, we proceed under the assumption that $L'_n = 0$.

We now let $L_{n+1} = \Sigma M_{[n+1,n+2)}$. Then $\alpha_n : \Sigma M_{\mathsf{I}} \to L_{n+1}$, $\beta_n : L_{n+1} \to \Sigma M_{\mathsf{I}}$, and we also have the inclusion $j : L_{n+1} \to \Sigma M_{\mathsf{I}}$ and the projection $q : \Sigma M_{\mathsf{I}} \to L_{n+1}$. The following product of six $2^{c-n}$-elementary automorphisms in the module $\Sigma M_{\mathsf{I}} \oplus L_{n+1}$

$$\begin{pmatrix} 1 & 0 \\ -q & 1 \end{pmatrix} \begin{pmatrix} 1 & j \\ 0 & 1 \end{pmatrix} \begin{pmatrix} 1 & 0 \\ -q & 1 \end{pmatrix} \begin{pmatrix} 1 & 0 \\ \alpha_n & 1 \end{pmatrix} \begin{pmatrix} 1 & -\beta_n \\ 0 & 1 \end{pmatrix} \begin{pmatrix} 1 & 0 \\ \alpha_n & 1 \end{pmatrix}$$

is easily seen to conjugate[13] the diagonal matrix $\mathrm{diag}(\beta_n\alpha_n, 0)$ into $\mathrm{diag}(jq, 0)$ (up to $2^{c-n}$-homotopy). Now we get the elementaries $\gamma^i$, $i = 2, 3, \cdots 7$, by applying $i_{n+1*}$ to each of the elementaries above, using (??) to identify $i_{n+1*}\Sigma M_{\mathsf{I}}$ with $M_{\mathsf{I}}$, and extending the resulting elementaries to be the identity on $(M'_{\mathsf{I}})^\perp$. With the identification mentioned, $i_{n+1*}(jq) : M_{\mathsf{I}} \to M_{\mathsf{I}}$ is the standard projection onto $M_{[n+1,n+2)}$, so the resulting $\varphi' = \gamma^7\gamma^6 \cdots \gamma^2\varphi : M' \to M'$ conjugates $\pi_{[n+1,\infty)} : M' \to M'$ into itself (up to a $2^{c-n}$-homotopy supported on $M'_{\mathsf{I}}$), and the proof is complete. $\square$

We can now finish the proof of Theorem ?? as follows. Recall that we have assumed that $y \in K_1^{cc}$ has $\partial(y) = 0$. Using the lemma, we can first write a representative of $y$ as $\varphi' \oplus \varphi : M' \oplus M \to M' \oplus M$ with $M' = M'_{[0,4)}$ and $M = M_{[4,\infty)}$. A standard Eilenberg swindle on the ray $(-\infty, 0]$ shows that $[cls\varphi'] = 0 \in K_1^{cc}$. Repeated applications of the lemma let us split $\varphi$ as an infinite direct sum

$$\varphi = \bigoplus_{k=1}^{\infty} \varphi_{4k} : \bigoplus M_{[4k,4k+4)} \to \bigoplus M_{[4k,4k+4)}$$

Indeed, the resulting (seven) infinite families of elementary automorphisms assemble to seven $cc$ elementary automorphisms on the direct sum, because of the conditions in Lemma ??, and a similar remark holds for the infinite family of stabilizations, and the infinite family of homotopies. The family $\{\Sigma\varphi_{4k}\}_k$ defines an element $z \in \lim^1 K_1$ and clearly $i(z) = y$.

## 6. COMPUTATIONS FOR TAME PAIRS.

Let $(\overline{Y}, C)$ be a compact, metric pair with $C$ tame in $\overline{Y}$ and with $\chi = (Y \xleftarrow{q} H \xrightarrow{p} C)$ the corresponding holink diagram. The groups $K_*(Y)$ and $K_*^{cc}(\overline{Y}, C; id_Y)$, as well

---

[13]This is based on the following (well known) fact: If, in some additive category, we have morphisms $\xi : X \to Y$ and $\eta : Y \to X$ satisfying $\xi\eta\xi = \xi$ and $\eta\xi\eta = \eta$, then the product

$$\begin{pmatrix} 1 & 0 \\ \xi & 1 \end{pmatrix} \begin{pmatrix} 1 & -\eta \\ 0 & 1 \end{pmatrix} \begin{pmatrix} 1 & 0 \\ \xi & 1 \end{pmatrix}$$

conjugates $diag(\eta\xi, 0)$ into $diag(0, \xi\eta)$.

as $K_{2-j}(C;p)_c$ and $\lim^1 K_{2-j}(C;p)_\epsilon$ for each $j \in \mathbb{N}$, are functorial in $(\overline{Y}, C) \in \mathcal{CM}^2_{tame}$, and we shall fit them into a natural chain complex

$$(69) \qquad K_2(Y) \xrightarrow{i_2} K_2^{cc}(\overline{Y}, C; id_Y) \xrightarrow{\Delta_2} K_1(C;p)_c \xrightarrow{\ell_1} K_1(Y) \xrightarrow{i_1} \cdots$$

and prove

**Theorem 6.1.** *Let* $j \in \mathbb{N}$. *The above chain complex is exact at* $K_{1-j}(C;p)_c$ *and* $K_{1-j}(Y)$, *and the homology group at* $K_{2-j}^{cc}(\overline{Y}, C; id_Y)$ *is naturally isomorphic to* $\lim^1 K_{2-j}(C;p)_\epsilon$.

**Remarks: 1)** This is the unreduced version of Theorem 1.2. All the ingredients in the following proof exist in reduced versions, so the proof of the other version is the same, mutatis mutandis.

**2)** We note that the theorem as stated contains no claims about the homology at $K_2^{cc}(\overline{Y}, C; id_Y)$ or exactness at $K_1(Y)$ and $K_1(C;p)_c$. In fact, two of the present authors (FXC and HJM) have worked out a definition of $\epsilon$-controlled groups $K_2(C;p)_\epsilon$ which appear to fit into an exact sequence

$$0 \to \lim^1 K_2(C;p)_\epsilon \to K_2^{cc}(\chi_+) \to K_1(C;p)_c \to 0.$$

Assuming this, the following proof will also cover the three cases currently excluded.
**3)** If we introduce "bumpy groups" $BK_*^{cc}$ (cf. [**?**]) by letting $BK_*^{cc}(\overline{Y}, C; id_Y) = \text{Ker}\,\Delta_*$, then the theorem translates into two exact sequences for each $j \in \mathbb{N}$, viz.

$$0 \to BK_{2-j}^{cc}(\overline{Y}, C; id_Y) \to K_{2-j}^{cc}(\overline{Y}, C; id_Y) \to K_{1-j}(C;p)_c \to K_{1-j}(Y),$$

and

$$K_{1-j}(C;p)_c \to K_{1-j}(Y) \to BK_{1-j}^{cc}(\overline{Y}, C; id_Y) \to \lim^1 K_{1-j}(C;p)_\epsilon \to 0.$$

They are analogous to the sequences in Siebenmann's Theorems I and II in [**?**]. In fact, if $Y$ is locally compact, metric with a metric one point compactification $\overline{Y} = Y \cup \{\infty\}$, and if $C = \{\infty\}$ is tame in $\overline{Y}$, then the connection to Siebenmann's exact sequences can be made more than an analogy. We shall not pursue this here, in part because we need the extension mentioned in **2)** in order to recover Siebenmann's Theorem II in the correct degree.
**4)** We recall that Chapman, see Theorems 13.1 and 13.2 of [**?**], studies the controlled boundary problem for a manifold $M^m$ parametrized over a compact metric space $B$ via a map $p : M \to B$. Assuming $m \geq 6$, and suitable "movability" and "tameness" conditions, he obtains a solution in terms of a primary obstruction $\sigma_\infty^p(M) \in \lim K_0(M - C_i)_{1/i}$ and a secondary obstruction $\tau_\infty^p(M) \in \lim^1 Wh(M - C_i)_{1/i}$. Here $\{C_i | i = 1, 2, 3 \cdots\}$ is an expanding and exhaustive sequence of compacta in $M$. For any space $Y$ with a map $p : Y \to B$, the group $Wh^p(Y)_\epsilon$ is defined geometrically in terms of $\epsilon$-controlled (in $B$) strong deformation retraction pairs, and the group $K_0^p(Y)_\epsilon$ is given as the transfer invariant subgroup of $Wh^{p\circ proj}(Y \times S^1)_\epsilon$. We may choose the sequence $C_i$ to be of the form $C_i = \rho^{-1}([0, 1 - 2^{-i}])$ for some proper map $\rho : M \to [0, 1)$. Now $\rho$ and $p$ together define a map $\bar{p} : M \to \overset{\circ}{c}(B)$, and we let $\overline{M} = M \amalg B$ with the teardrop topology. Urysohn's Metrizability Theorem can be used to show that $\overline{M}$ is metrizable, and presumaby Chapman's movability and tameness conditions imply that $B$ is tame in $\overline{M}$. Assuming this and the extension mentioned in **2)**, it is easy to produce natural maps

$$(70) \qquad \lim K_0(M - C_i)_{1/i} \to \tilde{K}_0^{cc}(B; \hat{p})_c,$$

and

$$(71) \qquad \lim{}^1 Wh(M - C_i)_{1/i} \to \lim{}^1 Wh(B; \hat{p})_\epsilon,$$

where $\chi(\overline{M}, B) = (M \xleftarrow{q} H \xrightarrow{\hat{p}} B)$. We conjecture that these maps are isomorphisms in some reasonable generality. If this is true, then it should be possible to give a direct definiton of a total obstruction $\varsigma(M; p) \in Wh(\overline{M}, B; id_M)$ such that $\Delta_1(\varsigma(M; p))$ correponds to $\sigma^p_\infty(M)$ and such that the homology class represented by $\varsigma(M; p)$ corresponds to $\tau^p_\infty(M)$ when $\sigma^p_\infty(M)$ vanishes. We plan to return to these questions in a later publication. $\qquad \square$

**Proof:** By Proposition **??** and Definition **??**, we have a long exact sequence

$$\cdots \longrightarrow K_{2-j}(Y) \xrightarrow{i_{2-j}} K_{2-j}^{cc}(\overline{Y}, C; id_Y) \xrightarrow{\bar{\Delta}_{2-j}} K_{2-j}^{cc}(\chi_+) \xrightarrow{\delta_{2-j}} K_{1-j}(Y) \longrightarrow \cdots.$$

It is valid for all $j \in \mathbb{Z}$, but right now we concentrate on $j \in \mathbb{N}$. Also, by Theorem **??**, the groups $K_{2-j}^{cc}(\chi_+)$ fit into short exact sequences

$$0 \to \lim{}^1 K_{2-j}(C; p)_\epsilon \xrightarrow{\gamma'_{2-j}} K_{2-j}^{cc}(\chi_+) \xrightarrow{\gamma''_{2-j}} K_{1-j}(C; p)_c \to 0$$

Since $\lim{}^1 K_{2-j}(\{c\}; c_Y)_\epsilon = 0$, a naturality argument applied to the morphism

$$(c_C, q, id_Y) : \chi \longrightarrow (\{c\} \xleftarrow{c_Y} Y \xrightarrow{id_Y} Y)$$

shows that the composition $\delta_{2-j} \circ \gamma'_{2-j} = 0$. Hence $\delta_{2-j}$ induces a homomorphism $\ell_{1-j} : K_{1-j}(C; p)_c \to K_{1-j}(Y)$, and if we take $\Delta_{2-j} = \gamma''_{2-j} \circ \bar{\Delta}_{2-j}$, an easy diagram chase establishes the desired chain complex (and the properties claimed), except for the first three terms. To get at these, on needs to show that the homomorphism $\bar{\Delta}_2 : K_2(\chi_+) \to K_1(Y)$ can be factored via $K_1(C; p)_c$. We leave this easy part of the extension mentioned in **2)** above, to the reader.

## 7. Ferry and Pedersen's vanishing result for $\lim{}^1 K_{2-j}(B; p)_\epsilon$.

Let $K$ be a finite, simplicial complex and let $\mathcal{A}$ be an additive category. Some time around 1987, Erik Kjær Pedersen, in a collaboration with Steve Ferry, gave a nice proof of the following result.

**Theorem 7.1.** *(Ferry and Pedersen - unpublished) There exists an $\epsilon > 0$ such that any element of $K_1(|K|; \mathcal{A})$ which can be represented by an automorphism of bound $\epsilon$ (of some geometric $\mathcal{A}$-object on $|K|$), can be "squeezed", i.e. it can be similarly represented with an arbitrarily small bound $\delta > 0$.*

Unfortunately, the proof is still unpublished, although there does exist an incomplete preprint, [**?**].

When Ferry and Pedersen first presented the proof at Odense University, two of the present authors (DRA and HJM) went through it in a variable coefficient version and found that it carries over with little extra trouble. Since the main source is unpublished, here we prefer to present the result as a conjecture, although we do feel it is a very strong such.

**Conjecture 7.2.** *Let $|K|$ be a finite polyhedron, and $p : E \to B = |K|$ be a map which has an iterated mapping cylinder stucture w.r.t. to the triangulation $K$ of $B$. Then the inverse systems $\{K_1(B; p)_\epsilon\}_\epsilon$ and $\{Wh(B; p)_\epsilon\}_\epsilon$ satisfy the Mittag Leffler condition. Hence the groups $\lim^1 K_{2-j}(B; p)_\epsilon$ and $\lim^1 \tilde{K}_{2-j}(B; p)_\epsilon$ vanish for each $j \in \mathbb{N}$.*

We hope that a final version of [**?**] will eventually be published.

Douglas R. Anderson
Dept. of Math., Syracuse University, Syracuse, NY 13210, USA
*E-mail address*: anderson@euclid.syr.edu

Francis X. Connolly
Dept. of Math., University of Notre Dame, Notre Dame, IN 46556, USA
*E-mail address*: francis.x.connolly.1.@nd.edu

Hans J. Munkholm
IMADA, Universitet, DK 5230 Odense M, Denmark
*E-mail address*: hjm@imada.ou.dk